\documentclass[11pt]{amsart}
\usepackage{mathrsfs}
\usepackage{amssymb}

\usepackage{titletoc}
\usepackage{amscd}
\usepackage{amsmath, amssymb}
\usepackage{amsfonts}
\usepackage[colorlinks,linkcolor=blue,citecolor=red, pdfstartview=FitH]{hyperref}

\numberwithin{equation}{section}

\theoremstyle{plain}
\newtheorem{thm}{Theorem}[section]
\newtheorem{theorem}[thm]{Theorem}
\newtheorem{lemma}[thm]{Lemma}
\newtheorem{corollary}[thm]{Corollary}
\newtheorem{proposition}[thm]{Proposition}

\newtheorem{conjecture}[thm]{Conjecture}
\theoremstyle{definition}

\newtheorem{remark}[thm]{Remark}

\newtheorem{definition}[thm]{Definition}

\newtheorem{example}[thm]{Example}

\newtheorem{defn-thm}[thm]{Definition-Theorem}

\newcommand{\sO}{{\mathcal O}}

\newcommand{\C}{{\mathbb C}}

\renewcommand{\P}{{\mathbb P}}

\newcommand{\R}{{\mathbb R}}

\renewcommand{\S}{{\mathbb S}}

\newcommand{\qtq}[1]{\quad\mbox{#1}\quad}
\newcommand{\bp}{\bar{\partial}}

\newcommand{\ts}{\otimes}

\newcommand{\btheorem}{\begin{theorem}}
\newcommand{\etheorem}{\end{theorem}}
\newcommand{\bproposition}{\begin{proposition}}
\newcommand{\eproposition}{\end{proposition}}
\newcommand{\bdefinition}{\begin{definition}}
\newcommand{\edefinition}{\end{definition}}
\newcommand{\bcorollary}{\begin{corollary}}
\newcommand{\ecorollary}{\end{corollary}}
\newcommand{\bproof}{\begin{proof}}
\newcommand{\eproof}{\end{proof}}
\newcommand{\bremark}{\begin{remark}}
\newcommand{\eremark}{\end{remark}}
\newcommand{\eexample}{\end{example}}
\newcommand{\bexample}{\begin{example}}

\newcommand{\elemma}{\end{lemma}}
\newcommand{\blemma}{\begin{lemma}}

\newcommand{\sq}{\sqrt{-1}}

\newcommand{\p}{\partial}

\renewcommand{\bar}{\overline}
\newcommand{\eps}{\varepsilon}

\renewcommand{\phi}{\varphi}

\newcommand{\ee}{\end{eqnarray*}}
\newcommand{\be}{\begin{eqnarray*}}

\newcommand{\beq}{\begin{equation}}
\newcommand{\eeq}{\end{equation}}

\newcommand{\bd}{\begin{enumerate}}
\newcommand{\ed}{\end{enumerate}}

\renewcommand{\hat}{\widehat}
\renewcommand{\tilde}{\widetilde}

\renewcommand{\>}{\rightarrow}

\newcommand{\bbe}{\bar{\beta}}

\newcommand{\al}{{\alpha}}
\newcommand{\ga}{{\gamma}}
\newcommand{\Ga}{{\Gamma}}

\usepackage{fancyhdr}
\renewcommand{\leftmark}{{\tiny RC-positivity and the  generalized energy density I: Rigidity}}

\begin{document}
\title{RC-positivity and the generalized energy density I: Rigidity} \makeatletter
\let\uppercasenonmath\@gobble
\let\MakeUppercase\relax
\let\scshape\relax
\makeatother

\author{Xiaokui Yang}
\date{}
\address{{Address of Xiaokui Yang: Morningside Center of Mathematics, Institute of
        Mathematics, Hua Loo-Keng Center of Mathematical Sciences,
        Academy of Mathematics and Systems Science,
        Chinese Academy of Sciences, Beijing, 100190, China.}}
\email{\href{mailto:xkyang@amss.ac.cn}{{xkyang@amss.ac.cn}}}

\noindent\thanks{\noindent This work was partially supported   by
China's Recruitment
 Program of Global Experts and  NSFC 11688101.
 }
\maketitle

\begin{abstract}  In this  paper,  we introduce a new energy density
function $\mathscr Y$ on the projective bundle $\P(T_M)\>M$ for a
smooth map $f:(M,h)\>(N,g)$ between Riemannian manifolds
$$\mathscr Y=g_{ij}f^i_\alpha f^j_\beta \frac{W^\alpha W^\beta}{\sum h_{\gamma\delta} W^\gamma W^\delta}.$$
  We get new Hessian
estimates to this energy density and obtain various new Liouville
type theorems for holomorphic maps, harmonic maps and pluri-harmonic
maps. For instance, we show that there is no non-constant
holomorphic map from a compact \emph{Hermitian manifold} with
positive (resp. non-negative) holomorphic sectional curvature to a
\emph{Hermitian manifold} with non-positive (resp. negative)
holomorphic sectional curvature.

\end{abstract}

{\small\setcounter{tocdepth}{1} \tableofcontents}

\section{Introduction}

Let $f:(M,h)\>(N,g)$ be a smooth map between two Riemannian
manifolds. In local coordinates $\{y^\alpha\}$ and $\{x^i\}$ on $M$
and $N$ respectively, there is an energy density function $e$ on $M$
$$e=|df|^2=g_{ij} h^{\alpha\beta}f^i_\alpha f^j_\beta. $$
Many milestone works are achieved in the last century by using
various  techniques in differential geometry and function theory in
analysis, and  thousands of mathematicians contributed significantly
in this rich field. There is a huge literature on the subject, and
we refer to the classical works \cite{Boc55, ES64, Yau75, Yau78,
Siu80, Yau82, JY93, MSY93, EL78, EL83, EL88, Xin96} and the
references therein.

  In this paper, we introduce a new energy density
function $\mathscr Y$ on the projective bundle $\P(T_M)\>M$ for the
smooth map $f:(M,h)\>(N,g)$ \beq \mathscr Y=g_{ij}f^i_\alpha
f^j_\beta \frac{W^\alpha W^\beta}{\sum h_{\gamma\delta} W^\gamma
W^\delta},\eeq which is motivated by
 the Leray-Grothendieck spectral sequence for abstract vector bundles used in our
previous paper \cite{Yang18a}. We obtain several new Hessian
estimates on this energy density (e.g. formulas (\ref{S1}),
(\ref{S2}), (\ref{S3}), (\ref{S11})). The key new ingredient is that
these new Hessian estimates can work for manifolds with partially
positive curvature tensors (e.g. the holomorphic sectional
curvature, or more generally, the RC-positivity for abstract vector
bundles introduced in \cite{Yang18}). In this paper, we only deal
with applications when $M$ is compact.

\subsection*{ Part I. The generalized energy density and rigidity of holomorphic maps} Let $f:(M,h)\>(N,g)$ be a holomorphic map between two
Hermitian
 manifolds. Let $\{z^\alpha\}_{\alpha=1}^m$ and $\{ \eta^i \}_{i=1}^n$ be the local holomorphic coordinates around $p\in M$ and $q=f(p)\in N$ respectively.
 We consider the generalized
 energy density
 \beq  \mathscr Y=g_{i\bar j}f_\alpha^i\bar{f}_{\beta}^j\frac{W^\alpha\bar W^\beta}{\sum h_{\gamma\bar\delta}W^\gamma \bar W^{\delta}} \label{density}\eeq
over the projective bundle $\P(T_M)\>M$ where $\{W^1,\cdots, W^m\}$
are the holomorphic coordinates on the fiber $T_pM$ with respect to
the given trivialization.  It is easy to see that $\mathscr Y$ is a
well-defined function on $\P(T_M)$. For simplicity, we set $\mathscr
H=\sum h_{\gamma\bar\delta}W^\gamma \bar W^{\delta}$. It is
well-known that $\mathscr H^{-1}$ is a  Hermitian metric on the
tautological line bundle $\sO_{T^*_M}(1)$ of the projective bundle
$\P(T_M)\>M$ (\cite[Lemma~9.1]{Gri65}). The complex Hessian of the
new energy density has the following estimate.

\btheorem\label{main0}Let $f:(M,h)\>(N,g)$ be a holomorphic map
between two Hermitian
 manifolds. We have the following inequality on the projective bundle
 $\P(T_M)\>M$,
\beq \sq\p\bp \mathscr  Y\geq \left(\sq \p\bp\log \mathscr
H^{-1}\right)\cdot \mathscr Y-\frac{\sq R_{i\bar j
k\bar\ell}f^i_{\alpha}\bar{f^j_{\beta}}f_\mu^k\bar{f_{\nu}^\ell}W^\mu\bar
W^\nu dz^\alpha\wedge d\bar z^\beta}{\mathscr H}. \label{S1} \eeq In
particular, if $f:M\>\C$ is a holomorphic function, then  \beq
\sq\p\bp \mathscr Y\geq \left(\sq \p\bp\log \mathscr
H^{-1}\right)\cdot \mathscr Y. \label{S-1}\eeq  \etheorem

\noindent As  applications of Theorem \ref{main0}, we obtain several
new rigidity theorems.

\btheorem\label{main} Let   $f:(M,h)\>(N,g)$  be a holomorphic map
between two Hermitian manifolds. Suppose $M$ is compact. If \bd
\item $(M,h)$ has positive (resp. non-negative) holomorphic sectional
curvature;

\item $(N,g)$ has non-positive (resp. negative) holomorphic sectional curvature,

\ed then $f$ is a constant map.

 \etheorem

\noindent Let's recall some classical results on  rigidity of
holomorphic maps.  The classical Schwarz-Pick lemma states that any
holomorphic map from the unit disc in the complex plane into itself
decreases the Poincar\'e metric. This was extended by Ahlfors
(\cite{Ahl38}) to maps from the disc into a hyperbolic Riemann
surface, and by Chern \cite{Che68} and Lu \cite{Lu68} to
higher-dimensional manifolds. A major advance was Yau's Schwarz
Lemma \cite{Yau78}, which says that a holomorphic map from a
complete K\"ahler manifold with \emph{Ricci curvature} bounded below
into a Hermitian manifold with holomorphic bisectional curvature
bounded above by a negative constant, is distance decreasing up to a
constant depending only on these bounds. In particular, there is no
nontrivial holomorphic map from compact K\"ahler manifolds with
positive Ricci curvature to
 Hermitian manifolds with non-positive holomorphic
bisectional curvature. By using the well-known ``Royden's trick" for
K\"ahler metrics, Royden was able to improve Yau's result  and show
that there is no nontrivial holomorphic map from compact K\"ahler
manifolds with positive Ricci curvature to
 \emph{K\"ahler manifolds} with non-positive holomorphic
sectional curvature (\cite{Roy80}).
 Later generalizations were mainly in two
directions: relaxing the curvature condition or the K\"ahler
assumption. A general philosophy is that  holomorphic maps from
``positively curved" complex manifolds to ``non-positively curved"
complex manifolds should be constant. \\

We confirmed in \cite[Theorem~1.7]{Yang18} a well-known problem of
Yau (\cite[Problem~47]{Yau82}) that if a compact K\"ahler manifold
$M$ has positive holomorphic sectional curvature, then it is
projective and rationally connected, i.e. any two points in $M$ can
be connected by a rational curve. On the other hand, if $N$ is Brody
hyperbolic, it has no rational curves. Hence, there is no
non-constant holomorphic maps from  compact K\"ahler manifolds $M$
with positive holomorphic sectional curvature to Brody hyperbolic
manifolds (e.g. Hermitian manifold $N$ with non-positive holomorphic
sectional curvature). Recently, Lei Ni also obtained a rigidity
theorem in \cite{Ni18} when $M$ and $N$ are both complete
\emph{K\"ahler manifolds} and one of the key ingredients is the
Royden's trick (\cite{Roy80}) for K\"ahler metrics. As it is shown
in Lemma \ref{main0} and Theorem \ref{main}, the new method in this
paper has the key advantage that they can work for Hermitian metrics
on both $M$ and $N$, and as it is well-known, the Royden's trick
does not work
on such manifolds.\\

  Theorem \ref{main} can \emph{not} be proved by using  algebraic methods developed in
\cite{Yang18a}. Indeed, when $(N,g)$ has
 non-positive holomorphic sectional curvature, the pullback  bundle $f^*(T^*_N)$
has no desired positivity as an abstract bundle.   On the other
hand, if $(N,g)$ has negative holomorphic sectional curvature, then
$(T_N, g)$ is RC-negative. However, it is obvious that the
RC-negativity is not preserved under the pull-back of $f$, unless
$f$ is a holomorphic submersion. As it is shown in
\cite[Corollary~3.8]{Yang18}, K\"ahler manifolds with negative first
Chern classes  have RC-negative tangent bundles and some of them can
contain rational curves, for instance, the quintic surface in
$\P^3$.

 As another application of the Hessian estimate in Theorem \ref{main0}, we
 obtain:

\bcorollary\label{main1}  Let $f:M\>N$ be a holomorphic map from a
compact complex manifold $M$ to a complex
 manifold $N$. If $\sO_{T^*_M}(-1)$ is RC-positive, and $N$ has non-positive
 holomorphic bisectional curvature, then $f$ is a constant map.
\ecorollary

\noindent  Let's explain the curvature condition on $M$ briefly. A
holomorphic line bundle $\mathscr L$ over a complex manifold $X$ is
called \emph{RC-positive} if it admits a smooth Hermitian metric
$h^{\mathscr L}$ such that its Chern curvature tensor $-\sq\p\bp\log
h^{\mathscr L}$ has at least one positive eigenvalue everywhere. The
RC-positivity of $\sO_{T_M^*}(-1)$ is a very weak curvature
condition. Indeed,
 we proved in \cite[Theorem~1.4]{Yang17}
that $\sO_{T^*_M}(-1)$ is RC-positive if and only if
$\sO_{T^*_M}(1)$ is not pseudo-effective. When restricted to the
fibers of $\P(T_M)$, one can deduce $\sO_{T^*_M}(-1)|_F\cong
\sO_{\P^{m-1}}(-1)$ is negative. \emph{Roughly speaking}, the
RC-positivity of $\sO_{T^*_M}(-1)$ over the projective bundle
$\P(T_M)$ means that $\sO_{T^*_M}(-1)$ has at least one ``positive
direction" along the base $M$ directions. Moreover, complex
manifolds $M$ with RC-positive tangent bundles have RC-positive
$\sO_{T^*_M}(-1)$ (\cite{Yang18a}). Recall that, a Hermitian
holomorphic vector bundle $(\mathscr E,h^{\mathscr E})$ over a
complex manifold $X$ is called \emph{RC-positive}, if for any $q\in
X$ and any nonzero vector $v\in \mathscr E_q$, there exists
\textbf{some} nonzero  vector $u\in T_qX$ such that \beq R^\mathscr
E(u,\bar u,v,\bar v)>0 .\eeq

\noindent There are many K\"ahler and non-K\"ahler complex manifolds
with RC-positive tangent bundles and we list some of them for
readers' convenience and for more details, we refer to \cite{Yang18,
Yang18a, Yang18b} and the references therein. \bd

\item[$\bullet$]  complex manifolds with positive holomorphic sectional
curvature;

\item[$\bullet$] Fano manifolds \cite[Corollary~3.8]{Yang18};
\item[$\bullet$]   manifolds with positive second Chern-Ricci
curvature \cite[Corollary~3.7]{Yang18};
\item[$\bullet$]  Hopf manifolds $\S^1\times \S^{2n+1}$ (\cite[formula (6.4)]{LY14});
\item[$\bullet$] products of complex manifolds with  RC-positive
tangent bundles.
 \ed
 We need to point out that Corollary
\ref{main1} can also be deduced from \cite{Yang18a} (see also
Corollary \ref{main2}).

\vskip 1\baselineskip

 We can define some other energy densities for  $f:(M,h)\>(N,g)$.   For instance,  on the projective bundle
 $\P(f^*T^*_N)\stackrel{\pi_1}{\>}M$, we
 have the energy density
\beq  \mathscr
Y_1=h^{\alpha\bar\beta}f_\alpha^i\bar{f}_{\beta}^j\frac{X_i\bar
X_j}{\sum g^{k\bar \ell} X_k \bar X_{\ell}} \label{density2}\eeq
 where
$\{X_i\}_{i=1}^n$ are the holomorphic coordinates on the fiber
$f^*T^*_qN$. For simplicity, we set $\mathscr H_1=\sum g^{k\bar
\ell} X_k \bar X_{\ell}$.

\bproposition\label{A1} Let  $f:(M,h)\>(N,g)$ be a holomorphic map
between two Hermitian manifolds. Then we have the following estimate
over $\P(f^*T^*_N)\stackrel{\pi_1}{\>}M$  \beq \sq\p\bp \mathscr
Y_1\geq \left(\sq \p\bp\log \mathscr H_1^{-1}\right)\cdot \mathscr
Y_1+\frac{\sq R_{\alpha\bar\beta \gamma\bar\delta}h^{\gamma\bar\nu}
h^{\mu\bar\delta }f_\mu^k\bar{f_{\nu}^\ell}X_k\bar X_\ell
dz^\alpha\wedge d\bar z^\beta}{\mathscr H_1}. \label{S2}\eeq
\eproposition

\noindent Similarly, we can define \beq \mathscr
Y_2=f_\alpha^i\bar{f}_{\beta}^j\frac{X_i\bar X_j}{\sum g^{k\bar
\ell} X_k \bar X_{\ell}} \cdot \frac{W^\alpha\bar W^\beta}{\sum
h_{\gamma\bar\delta}W^\gamma \bar W^{\delta}} \label{density3}\eeq
over the projective bundle $\P(\pi^*f^*T^*_N)\>\P(T_M)$ where
$\pi:\P(T_M)\>M$ is the natural projection. Recall that $\mathscr
H=\sum h_{\gamma\bar\delta}W^\gamma \bar W^{\delta}$ and $\mathscr
H_1=\sum g^{k\bar \ell} X_k \bar X_{\ell}$.

\bproposition\label{A2} Let  $f:(M,h)\>(N,g)$ be a holomorphic map
between two Hermitian manifolds. Then we have the following
inequality over $\P(\pi^*f^*T^*_N)\>\P(T_M)$   \beq \sq\p\bp
\mathscr Y_2\geq \left(\sq \p\bp\log \mathscr H^{-1}+\sq \p\bp\log
\mathscr H_1^{-1}\right)\cdot \mathscr Y_2. \label{S3}\eeq
\eproposition

\vskip 1\baselineskip

\noindent As applications of \emph{generalized versions} of Theorem
\ref{main0}, Proposition \ref{A1} or Proposition \ref{A2}, and some
deep approximations established by
 Demailly-Peternell-Schneider (\cite{DPS94}), we conclude:

\bcorollary\label{main2} Let $f:M\>N$ be a holomorphic map between
compact complex
 manifolds. If $\sO_{T^*_M}(-1)$ is RC-positive and $\sO_{T^*_N}(1)$ is nef, then $f$ is a constant map.
\ecorollary

\noindent Note that $\sq \p\bp\log \mathscr H_1^{-1}$ is the
curvature tensor of the line bundle
$\left(\sO_{f^*(T_N)}(-1),\mathscr H_1\right)$.  When $\dim_\C N\geq
2$,  $f^*\left(\sO_{T^*_N}(1)\right)$ and $\sO_{f^*T_N}(-1)$ are not
isomorphic. Actually, the restriction of $\sO_{f^*(T_N)}(-1)$ on
each fiber $\P(f^*T^*_qN)\cong \P^{n-1}$ is isomorphic to
$\sO_{\P^{n-1}}(-1)$. Hence, $\sO_{f^*(T_N)}(-1)$ can not be nef in
any case. It is easy to see that if $(N,h)$ has non-positive
holomorphic bisectional curvature, then $\sO_{f^*(T_N)}(-1)$ is
RC-nonnegative along the base $M$ directions. As analogous to the
proof of Corollary \ref{main1}, it is not hard to see that we can
 prove Corollary \ref{main2} by using the generalized version of inequality
 (\ref{S1}) or (\ref{S3}) for \emph{some other new energy density}  since
 the desired asymptotic
 metrics
  constructed in
\cite[Theorem~1.12]{DPS94} are on the vector bundle
$\mathrm{Sym}^{\ts k}T^*_N$. It worths to point out that, Corollary
\ref{main2} was firstly established in \cite[Theorem~1.1]{Yang18a}
 by using the Leray-Grothendieck spectral
sequences and isomorphisms of various cohomology groups for
\emph{abstract vector bundles}. The algebraic proof in
\cite{Yang18a} is much more effective in
this setting!\\

  In the same spirit of Theorem \ref{main} and Corollary \ref{main2}, we
  have:

\bcorollary\label{main3} Let $f:M\>N$ be a holomorphic map between
complex
 manifolds. Suppose $M$ is compact. If \bd \item  $\sO_{T_M}(1)$ is nef;

 \item  $N$ has a Hermitian metric with negative holomorphic sectional curvature, \ed then $f$ is a constant map.
\ecorollary

\noindent Recently, there are many important works on holomorphic
sectional curvature with various positivity in the K\"ahler setting,
and we refer to \cite{HLW10, HW12,HW12, HLW14, HW15, ACH15, Liu16,
WY16, WY16a, Nom16, Yang16, YZ16, YbZ16, AHZ16, AH17,Yang18,
Yang18a, Yang18b, NZ18a, NZ18b, Ni18, Mat18a, Mat18b, Gue18, Zha18}
and the references therein.  The results in this paper demonstrate
certain similarity  between Hermitian metrics and K\"ahler metrics
with such curvature positivity. It might be an interesting problem
to ask whether all compact Hermitian manifolds with positive
holomorphic sectional curvature are K\"ahler or projective.

 \vskip 1\baselineskip

\subsection*{Part II. Rigidity of harmonic maps and pluri-harmonic maps.}

 Let $f:(M,h)\>(N,g)$ be a smooth map from
a Hermitian
 manifold to a Riemannian manifold. Let $\{z^\alpha\}_{\alpha=1}^m$ and $\{ x^i \}_{i=1}^n$ be the
 local holomorphic coordinates and real coordinates on $M$ and $N$ respectively. We consider the generalized
 energy density
 \beq  \mathscr Y=g_{i j}f_\alpha^i\bar{f}_{\beta}^j\frac{W^\alpha\bar W^\beta}{\sum h_{\gamma\bar\delta}W^\gamma \bar W^{\delta}} \eeq
over the projective bundle $\P(T_M)$, where $(W^1,\cdots, W^m)$ are
the holomorphic coordinates on the fiber $T_pM$. We set $\mathscr
H=h_{\gamma\bar\delta}W^\gamma \bar W^{\delta}.$\\

 It is well-known that
pluri-harmonic maps are generalizations of holomorphic maps and
harmonic maps. Indeed, a smooth map $f$ from a complex manifold $M$
to a K\"ahler (or Riemannian) manifold  $N$  is
\emph{pluri-harmonic} if and only if for any holomorphic curve $i$:
$C\>M$, the composition $f\circ i$ is harmonic. On the other hand,
$\pm$ holomorphic maps between K\"ahler manifolds are
pluri-harmonic. As analogous to Theorem \ref{main0}, we obtain:

\btheorem\label{main10} If $f:(M,h)\>(N,g)$ is a pluri-harmonic
harmonic map, then we have the following inequality on $\P(T_M)$
\beq \sq\p\bp \mathscr Y\geq \left(\sq \p\bp\log \mathscr
H^{-1}\right)\cdot \mathscr Y-\frac{\sq R_{i\ell k j}
f^i_{\alpha}{f^j_{\bar\beta}}f_\gamma^k{f_{\bar\delta}^\ell}W^\gamma\bar
W^\delta dz^\alpha\wedge d\bar z^\beta}{\mathscr H}. \label{S11}
\eeq \etheorem

As applications of Theorem \ref{main10}, we study the geometry of
compact K\"ahler manifolds with \emph{negative/non-positive
Riemannian sectional curvature} by using harmonic maps and
pluri-harmonic maps into such manifolds.  The motivation is the
following conjecture proposed by S.-T. Yau
(\cite[Problem~37]{Yau82}) :

\begin{conjecture}\label{Conj}  Let $(X,\omega)$ be a compact K\"ahler manifold with $\dim_\C
X>1$. Suppose $(X,\omega)$ has negative Riemannian sectional
curvature, then $X$ is rigid, i.e. $X$ has only one  complex
structure.
\end{conjecture}

\noindent It is a fundamental problem on the rigidity of K\"ahler
manifolds with negative curvature.  S.-T. Yau proved in
\cite[Theorem~6]{Yau77} that when $X$ is covered by a $2$-ball, then
any complex surface oriented homotopic to $X$ must be biholomorphic
to $X$. By using the terminology of ``strongly negativity", Y.-T.
Siu established in \cite[Theorem~2]{Siu80} that a compact
 K\"ahler manifold of the same homotopy type as  a compact K\"ahler manifold
$(X,\omega)$ with strongly negative curvature and $\dim_{\C}X>1$
must be either biholomorphic or conjugate biholomorphic to $X$.  It
is well-known that the strongly negative curvature condition can
imply the negativity of the Riemannian sectional curvature. Hence,
 Conjecture \ref{Conj} holds under that stronger curvature condition. Note
that when $\dim_{\C}X=2$, Conjecture \ref{Conj} has been completely
solved in \cite{Zhe95} by F.-Y. Zheng. When $\dim_\C X\geq 3$,
Conjecture \ref{Conj} is still widely open since there is no
effective method to deal with the Riemannian sectional curvature on
complex manifolds.\\

Before stating the applications of Theorem \ref{main10}, we recall
the strategy in establishing Siu's strong rigidity mentioned above.
Let $f:(M,h)\>(N,g)$ be a smooth map between two compact K\"ahler
manifolds.

\bd \item[(A)] Suppose $(N,g)$  has strongly non-positive curvature
in the sense of Siu. If $f$ is a harmonic map,  then $f$ is
pluri-harmonic (see Lemma \ref{pluri});

\item[(B)] Suppose $(N,g)$ has strongly negative curvature in the
sense of Siu.  Then any pluri-harmonic map $f:M\>(N,g)$ is
holomorphic or anti-holomorphic provided  $\mathrm{rank}_{\R} df\geq
4$ (see Remark \ref{siure}). \ed

\noindent As inspired by Yau's conjecture \ref{Conj} and Siu's
strategy in steps (A) and (B), we attempt to investigate Riemannian
(or K\"ahler) manifolds $(N,g)$ with non-positive Riemannian
sectional curvature. As an application of Theorem \ref{main10}, we
obtain:

\btheorem\label{main11}  Let $f:M\>(N,g)$ be a pluri-harmonic map
from a compact complex manifold $M$ to a Riemannian
 manifold $(N,g)$ with \emph{non-positive
 Riemannian sectional curvature}. If $\sO_{T^*_M}(-1)$ is RC-positive, then $f$ is a constant map.
\etheorem

\noindent  Theorem \ref{main11} still holds when the target
 manifold is K\"ahler:

\btheorem\label{main12} Let $f:M\>(N,g)$ be a pluri-harmonic map
from a compact complex manifold $M$ to a K\"ahler manifold $(N,g)$
with \emph{non-positive Riemannian sectional curvature}. If
$\sO_{T^*_M}(-1)$ is RC-positive, then $f$ is a constant map.
\etheorem

\noindent  As we pointed out before, there are many K\"ahler and
non-K\"ahler manifolds with RC-positive $\sO_{T^*_M}(-1)$. In
particular, we get: \bcorollary\label{coro} Let $f:M\>(N,g)$ be a
pluri-harmonic map from a compact complex manifold $M$ to a K\"ahler
or Riemannian
 manifold $(N,g)$ with \emph{non-positive
 Riemannian sectional curvature}. If $M$ has a Hermitian metric with positive holomorphic sectional curvature, then $f$ is a constant map.
\ecorollary

\noindent As an application of Theorem \ref{main11} and ideas in
\cite{Siu80, Sam85, Sam86, YZ91, JY91, JY93,
 LY14, WY18},  we consider harmonic maps from complex manifolds into
 Riemannian manifolds.

\btheorem\label{main13}  Let $f:(M,h)\>(N,g)$ be a harmonic map from
a compact
 K\"ahler manifold $(M,h)$ to a Riemannian
 manifold $(N,g)$ with \emph{non-positive
complex sectional curvature}. If $\sO_{T^*_M}(-1)$ is RC-positive,
then $f$ is a constant map. \etheorem

\noindent  Theorem \ref{main13} still holds if $f$ is a Hermitian
harmonic map from an astheno-K\"ahler manifold $(M,h)$ (i.e.
$\p\bp\omega^{m-2}_h=0$) (introduced by Jost-Yau in \cite{JY93}) to
a Riemannian
 manifold $(N,g)$ with non-positive
complex sectional curvature.  In \cite{WY18}, we obtained some
results by using the classical Chern-Lu type inequality
(\ref{hessian}) under a much stronger condition that $T_M$ is
uniformly RC-positive.

\bremark As analogous to the classical theory of harmonic maps,
there are many further applications of the generalized energy
density (\ref{density}). They are discussed briefly in Section
\ref{further}. For instances, \bd

\item[$\bullet$] the RC-positivity for Riemannian curvature tensor;

 \item[$\bullet$] the extension of Yau's  function
theory on complete manifolds \cite{Yau75, Yau78};

\item[$\bullet$] the first and second variations of the generalized  energy functions, and
the applications in investigating the existence of rational curves
on manifolds with RC-positive curvature by analytical methods
(\cite{SU81, SY80});

\item[$\bullet$] the analytical extension of  methods in this paper  to hyperbolic manifolds;

\item[$\bullet$] the generalized energy density on Grassmannian manifolds
$\mathrm{Gr}(k,T_M)$ for RC-positivity in $k$ linearly independent
directions. \ed
 The Ricci-flow and K\"ahler-Ricci flow approaches (\cite{Ham82,
Cao85})  in this setting and  the applications of parabolic
estimates corresponding to  formulas (\ref{S1}), (\ref{S2}),
(\ref{S3}) and (\ref{S11})  are also expectable. The details of some
topics listed above will appear somewhere else. \eremark

This paper is organized as follows. In Section \ref{2}, we describe
the relationship between the classical energy identity and the
generalized energy density.  In Section \ref{3},  we prove Theorem
\ref{main0}, Theorem \ref{main},
 Corollary \ref{main1} and Corollary \ref{main2}. In Section
 \ref{4}, we investigate harmonic maps and pluri-harmonic from complex manifolds to Riemannian manifolds and K\"ahler manifolds, and establish Theorem \ref{main10}, Theorem \ref{main11},
 Theorem \ref{main12} and Theorem \ref{main13}.\\

\noindent \textbf{Acknowledgements.} The author is very grateful to
Professor
 K.-F. Liu and Professor S.-T. Yau for their support, encouragement and stimulating
discussions over  years. The author  would also like to thank B.-L.
Chen, F.-Q. Fang, N. Mok, J. Wang, V. Tosatti,  W.-P. Zhang and
X.-Y. Zhou for some useful suggestions.

\vskip 2\baselineskip

\section{The classical energy density and the generalized energy
density}\label{2}

\subsection{Holomorphic maps between complex manifolds.}
 Let $f:(M,h)\>(N,g)$ be a smooth map between two
Hermitian
 manifolds. Let $\{z^\alpha\}_{\alpha=1}^m$ and $\{ \eta^i \}_{i=1}^n$ be the local holomorphic coordinates around $p\in M$ and $q=f(p)\in N$ respectively.
The classical $\p$-energy density is defined as \beq u=|\p
f|^2=g_{i\bar
j}h^{\alpha\bar\beta}f_\alpha^i\bar{f}_{\beta}^j.\label{S10}\eeq
Here $\p f$ can be regarded as a section of the complex vector
bundle $E=T^*_M\ts f^*T_N$ and $u$ is the norm square of $\p f$ with
respect to the induced metric on $E$.  As it is well-known, if $f$
is a holomorphic map, by using standard Bochner technique, one has
the following Chern-Lu inequality (\cite{Che68,Lu68}, see also
\cite[Lemma~5.1]{Yang18a}).
\begin{lemma}\label{Lemma4.1}
Let $f:(M,h)\rightarrow (N,g)$ be a holomorphic map between two
Hermitian manifolds. Then \beq  \sq \p\bp u\geq \sq \left(
R_{\alpha\bar\beta\gamma\bar\delta}h^{\mu\bar\delta}h^{\gamma\bar\nu}g_{i\bar
j}f^i_\mu\bar{f^j_\nu} -R_{i\bar j
k\bar\ell}f^i_{\alpha}\bar{f^j_{\beta}}\left(h^{\mu\bar\nu}f_\mu^k\bar{f_{\nu}^\ell}\right)\right)dz^\alpha\wedge
d\bar z^\beta,\label{S01} \eeq

\noindent and \beq \mathrm{tr}_{\omega_h}\left(\sq \p\bp
u\right)\geq \left(h^{\alpha\bar\beta}
R_{\alpha\bar\beta\gamma\bar\delta}\right)h^{\mu\bar\delta}h^{\gamma\bar\nu}\left(g_{i\bar
j}f^i_\mu\bar{f^j_\nu}\right) -R_{i\bar j
k\bar\ell}\left(h^{\alpha\bar\beta}f^i_{\alpha}\bar{f^j_{\beta}}\right)\left(h^{\mu\bar\nu}f_\mu^k\bar{f_{\nu}^\ell}\right).
\label{S02}\eeq
\end{lemma}

\noindent Formulas (\ref{S01}), (\ref{S02}) and their parabolic
analogs  have many fantastic applications in differential geometry
and the theory of Ricci flows, and we refer to \cite{Che68, Lu68,
Yau78, Roy80, EL83, EL88, Tos07, Che09, YZ16, Ni18} and the
references therein.

\bremark In the formula \ref{S01}, if we choose
$h_{\alpha\bar\beta}=\delta_{\alpha\bar\beta}$ at some point, then
the curvature term on the domain manifold is $\sq
R_{\alpha\bar\beta\gamma\bar\delta}g_{i\bar
j}f^i_\gamma\bar{f^j_\delta}dz^\alpha\wedge d\bar z^\beta$. In the
approach by using maximum principle,  it is hard to get the desired
positivity from this term. \eremark

Next, we introduce the generalized energy density on the projective
bundle $\P(T_M)$. The points of the projective bundle $\P(T_M)$ of
$T^*_M\>M$ can be identified with the hyperplanes of $T^*_M$. Note
that every hyperplane $ V$ in $T^*_pM$ corresponds bijectively to
the line of linear forms in $T_pM$ which vanish on $ V$.
 Let $\pi : \P(T_M) \> M$ be
the natural projection. Suppose $(W^1,\cdots, W^m)$ are the
holomorphic coordinates on the fiber of $TM$. The generalized
$\p$-energy density on the projective bundle $\P(T_M)\>M$ is defined
as
 \beq  \mathscr Y=g_{i\bar j}f_\alpha^i\bar{f}_{\beta}^j\frac{W^\alpha\bar W^\beta}{\sum h_{\gamma\bar\delta}W^\gamma \bar W^{\delta}}. \label{density01}\eeq
It is easy to easy see that (\ref{density01}) is well-defined on
$\P(T_M)$. The classical energy density (\ref{S10}) and the
generalized energy density (\ref{density01}) are related in the
following way.

\bproposition\label{relation} We have the following relation \beq
|\p f|^2=m\pi_*(\mathscr Y), \eeq where $\pi_*$ is the fiberwise
integration with respect to the fiberwise  Fubini-Study metric
$\omega=-\sq \p\bp\log (\sum h_{\gamma\bar\delta}W^\gamma \bar
W^{\delta})$.
 \eproposition \bproof Indeed, \be \pi_*(\mathscr
Y)=\int_{\P(T_pM)} g_{i\bar
j}f_\alpha^i\bar{f}_{\beta}^j\frac{W^\alpha\bar W^\beta}{\sum
h_{\gamma\bar\delta}W^\gamma \bar W^{\delta}}\cdot
\frac{\omega^{m-1} }{(m-1)!}=  g_{i\bar
j}f_\alpha^i\bar{f}_{\beta}^j\cdot \frac{
h^{\alpha\bar\beta}}{m}=\frac{|\p f|^2}{m},\ee where we use the
well-known identity for the Fubini-Study metric
$\omega_{\mathrm{FS}}$ on $\P^{m-1}$ \beq \int_{\P^{m-1}}
\frac{W^\alpha\bar W^\beta}{|W|^2}
\frac{\omega_{\mathrm{FS}}^{m-1}}{(m-1)!}=
\frac{\delta^{\alpha\bar\beta}}{m}. \eeq \eproof

\noindent
 We can also define a conformal change for the generalized energy density
 \beq \mathscr Y_\phi=e^\phi\mathscr Y =e^\phi g_{i\bar j}f_\alpha^i\bar{f}_{\beta}^j\frac{W^\alpha\bar W^\beta}{\sum h_{\gamma\bar\delta}W^\gamma \bar W^{\delta}} \label{density02}\eeq
for any $\phi\in C^\infty(\P(T_M),\R)$. As we pointed out before,
$\mathscr H=\sum h_{\gamma\bar\delta}W^\gamma \bar W^{\delta}$ is a
(local) Hermitian metric on $\sO_{T^*_M}(-1)$, and so $\mathscr H
e^{-\phi}$ is  also a Hermitian metric on  $\sO_{T^*_M}(-1)$.
Actually, any Hermitian metric on it takes such a form.

\subsection{Pluri-harmonic maps into
Riemannian manifolds} Let $(M.h)$ be a complex Hermitian manifold,
$(N,g)$ be a Riemannian manifold and $f:M\>N$ be a smooth map.  We
denote by $\mathscr E=f^*(TN)$ and endow it with the induced
{Levi-Civita connection} $\nabla^\mathscr E$ from $TN$. There is a
natural decomposition $\nabla^{\mathscr E}=\bp_\mathscr
E+\p_\mathscr E$ according to the complex structure on $M$. Let
$\{z^\alpha\}_{\alpha=1}^m$ be the local holomorphic coordinates on
$M$ and $\{x^i\}_{i=1}^n$ the local coordinates on $N$. Let
$e_i=f^*(\frac{\p}{\p x^i})$. There are three $\mathscr E$-valued
$1$-forms, i.e. \beq \bp f =\frac{\p f^i}{\p \bar z^\beta}d\bar
z^\beta\ts e_i, \ \ \p f=\frac{\p f^i}{\p z^\alpha}dz^\alpha\ts
e_i,\ \ \  df=\bp f+\p f. \eeq It is easy to see that \beq u=|\p
f|^2=g_{i j}h^{\alpha\bar\beta}f_\alpha^i\bar{f}_{\beta}^j,
\qtq{and} u=\frac{1}{2}|df|^2= |\bp f|^2=|\p f|^2 .\eeq $f$ is
called a \emph{harmonic map} if it is a critical point of the
Euler-Lagrange
equation of the total energy $E(f)=\int_M u dV_M$.\\

 We consider the generalized
 energy density
 \beq  \mathscr Y=g_{i j}f_\alpha^i\bar{f}_{\beta}^j\frac{W^\alpha\bar W^\beta}{\sum h_{\gamma\bar\delta}W^\gamma \bar W^{\delta}} \label{density03}\eeq
over the projective bundle $\P(T_M)$, where $(W^1,\cdots, W^m)$ are
the holomorphic coordinates on the fiber $T_pM$. We set $\mathscr
H=h_{\gamma\bar\delta}W^\gamma \bar W^{\delta}.$ Similarly, we have
\beq |df|^2=2m\pi_*(\mathscr Y).\eeq

\bdefinition  A smooth map $f: M\>(N,g)$ from a complex manifold to
a Riemannian manifold  is called \emph{pluri-harmonic} if it
satisfies $\p_\mathscr E\bp f=0$, i.e. \beq \left(\frac{\p^2 f^i}{\p
z^\alpha\p\bar z^\beta}+\Gamma_{jk}^i\frac{\p f^j}{\p \bar
z^\beta}\frac{\p f^k}{\p z^\alpha}\right)dz^\alpha\wedge d\bar
z^\beta\ts e_i=0\label{pluriharmonic}\eeq where $\Gamma_{jk}^i$ is
the Christoffel symbol of the Levi-Civita connection on $(N,g)$.
 \edefinition
\noindent It is easy to see that $\pm$ holomorphic maps from complex
manifolds to K\"ahler manifolds are pluri-harmonic.

 Let $(N,g)$ be a Riemannian manifold with the Levi-Civita connection
$\nabla$. Its  curvature tensor is defined as $$
R(X,Y)Z=\nabla_X\nabla_YZ-\nabla_Y\nabla_XZ-\nabla_{[X,Y]}Z$$ for
any $X,Y,Z\in \Gamma(N,TN)$. In the local coordinates $\{x^i\}$ of
$N$, we adopt the convention $$
R(X,Y,Z,W)=g(R(X,Y)Z,W)=R_{ijk\ell}X^iY^jZ^kW^\ell.$$ It is easy to
see that \beq R_{ijk}^\ell=\frac{\p\Gamma_{k j}^\ell}{\p
x^i}-\frac{\p\Gamma_{k i}^{\ell}}{\p x^j}+\Gamma^{p}_{k
j}\Gamma^{\ell}_{p i}-\Gamma^{p}_{k i}\Gamma^{\ell}_{p j}, \qtq{and}
R_{ijk\ell}=g_{s\ell}R_{ijk}^s. \label{rcurvature}\eeq

 The following constraint equation for pluri-harmonic maps is
well-known (e.g. \cite[Lemma~1.3]{OU91} and \cite[Lemma~2.4]{WY18}),
and it follows by taking first order derivatives on the
pluri-harmonic equation (\ref{pluriharmonic}). \blemma If
$f:M\>(N,h)$ is a pluri-harmonic map, then
 \beq
 R_{ikj\ell}f^i_\alpha f^j_{\bar\beta} f^k_\gamma=0
\label{D}
  \eeq
for any $\alpha,\beta,\gamma$ and $\ell$ where $f^i_\alpha=\frac{\p
f^i}{\p z^\alpha}$ and $f^j_{\bar \beta}=\frac{\p f^j}{\p\bar
z^\beta}$. \elemma

\bremark By using the constraint equation (\ref{D}), we have \beq
\hat C=R_{ik\ell
j}\left(h^{\alpha\bar\beta}f^i_{\alpha}f^j_{\bar\beta}\right)
\left(h^{\gamma\bar\delta}f^k_{\gamma}f^\ell_{\bar\delta}\right)=0.\label{hatC}
\eeq If $(N,g)$ has  positive or negative constant Riemannian
sectional curvature, one can deduce from (\ref{hatC}) that
$\mathrm{rank}_\R df\leq 2$ (\cite{Sam85, Sam86}). \eremark

\noindent As analogous to the Chern-Lu inequalities for holomorphic
maps in Lemma \ref{Lemma4.1},  Wang and the author obtained in
\cite[Proposition~3.2]{WY18} the following inequalities for
pluri-harmonic maps.

\blemma\label{cri} Let $f:(M,h)\>(N,g)$ be a pluri-harmonic map from
a
 Hermitian manifold $M$ to a Riemannian manifold $(N,g)$.
Then we have \beq \sq \p\bp u \geq  \sq\left( R_{\alpha\bar \beta
\gamma\bar \delta} h^{\gamma \bar \nu} h^{\mu\bar \delta} g_{ij}
f^{i }_\mu {f^{j}_{\bar\nu}} - R_{i\ell k j}f^i_{\alpha}f^j_{\bar
\beta}
\left(h^{\gamma\bar\delta}f^k_{\gamma}f^\ell_{\bar\delta}\right)\right)dz^\alpha\wedge
d\bar z^\beta.\label{hessian} \eeq

\noindent
 and \beq \mathrm{tr}_{\omega_h}\left(\sq \p\bp
u\right)\geq \left(h^{\alpha\bar\beta}
R_{\alpha\bar\beta\gamma\bar\delta}\right)h^{\mu\bar\delta}h^{\gamma\bar\nu}\left(g_{i
j}f^i_\mu{f^j_{\bar\nu}}\right) -R_{i\ell k
j}\left(h^{\alpha\bar\beta}f^i_{\alpha}{f^j_{\bar{\beta}}}\right)\left(h^{\mu\bar\nu}f_\mu^k{f_{\bar\nu}^\ell}\right).
\label{hessian2} \eeq \elemma

\bremark The formulations of curvature terms
 in Lemma \ref{cri} and Lemma
\ref{Lemma4.1} on the target manifolds are different. \eremark


For more detailed discussions on various harmonic maps and
pluri-harmonic maps, we refer to
 \cite{EL78, Siu82, EL83, Ohn87, EL88, Uda88, CT89, OU90, OV90, YZ91, MSY93, Uda94, JZ97, Lou99, Ni99, DEP03,
 Tos07, Zhang12,
  YHD13, Dong13, LY14, Sch15, YZ16, Yang18a, Yang18b, Ni18, WY18} and the references therein.

\vskip 2\baselineskip

\section{The Hessian estimates and rigidity of holomorphic maps}\label{3}
 In this section, we prove Theorem \ref{main0}, Theorem \ref{main},
 Corollary \ref{main1} and Corollary \ref{main2}. The proofs of Proposition \ref{A1} and Proposition \ref{A2} are similar to that in Theorem \ref{main0}. To avoid identifications in algebraic geometry,  we shall use
straightforward local
computations  for readers' convenience.\\

\noindent \emph{The proof of Theorem \ref{main0}.}   For simplicity,
we write $\mathscr F=\mathscr Y\cdot \mathscr H$. A straightforward
calculation on $\P(T_M)$ yields \beq \p\bp\mathscr Y=\left(\p\bp\log
\mathscr H^{-1}\right) \mathscr Y+\frac{\p\bp \mathscr F}{\mathscr
H}+ \frac{\mathscr F\p \mathscr H\wedge \bp \mathscr H}{\mathscr
H^3}-\frac{\p \mathscr H\wedge \bp \mathscr F+\p\mathscr F\wedge \bp
\mathscr H }{\mathscr H^2}.\label{E0}  \eeq  Moreover, we set
$F^i=f^i_\alpha
 W^\alpha$ and so $\mathscr F=g_{i\bar j}F^i\bar F^j$.  Since $f$ is holomorphic, we have
\beq\p \mathscr F=\p g_{i\bar j} \cdot F^i\bar{F^j}+g_{i\bar j}\cdot
\p F^i\cdot \bar F^j,\ \ \ \ \bp\mathscr F=\bp g_{i\bar j}\cdot
F^i\bar{F^j}+ g_{i\bar j} F^i\bar{\p F^j}. \eeq

\noindent  Therefore, we deduce
\begin{eqnarray} &&\label{E3}\frac{\p\bp \mathscr F}{\mathscr H}+ \frac{\mathscr F\p
\mathscr H\wedge \bp \mathscr H}{\mathscr H^3}-\frac{\p \mathscr
H\wedge \bp \mathscr F+\p\mathscr F\wedge \bp \mathscr H }{\mathscr
H^2}\\ \nonumber&=&\frac{\p\bp g_{i\bar j}\cdot F^i\bar{F^j}+ \p
g_{i\bar j}\cdot F^i\cdot \bar{\p F^j}-\bp g_{i\bar j}\cdot\p
F^i\cdot \bar{F^j}+ g_{i\bar j}\p F^i\bar{\p F^j}}{\mathscr H}\\
\nonumber&&+\frac{\mathscr F \p\log\mathscr H^{-1}\wedge
\bp\log\mathscr H^{-1}}{\mathscr H}\\
\nonumber&&-\frac{\p\log\mathscr H^{}\cdot \left(\bp g_{i\bar
j}\cdot F^i\bar{F^j}+ g_{i\bar j} F^i\bar{\p F^j}\right)+\left(\p
g_{i\bar j} \cdot F^i\bar{F^j}+g_{i\bar j}\cdot \p F^i\cdot
F^j\right)\cdot \bp\log\mathscr H^{}}{\mathscr H}.\end{eqnarray}

\noindent  We define a   $(1,1)$-form on $\P(T_M)$: \beq \mathscr
W=g_{i\bar j}\left(\p F^i+F^i\p\log\mathscr
H^{-1}+T^i\right)\wedge\bar{ \left(\p F^j+F^j\p\log\mathscr
H^{-1}+T^j\right)} \label{E2} \eeq where $T^i=g^{i\bar\ell}\frac{\p
g_{k\bar\ell}}{\p z^p}F^k\p f^p$. It is easy to see that \beq
g_{i\bar j} T^i=\frac{\p g_{k\bar j} }{\p z^p} F^k\p f^p=(\p
g_{k\bar j}) F^k.\label{E21}\eeq

\noindent By using (\ref{E3}), (\ref{E2}) and (\ref{E21}), we obtain
\begin{eqnarray}&&\nonumber\frac{\p\bp \mathscr F}{\mathscr H}+ \frac{\mathscr F\p
\mathscr H\wedge \bp \mathscr H}{\mathscr H^3}-\frac{\p \mathscr
H\wedge \bp \mathscr F+\p\mathscr F\wedge \bp \mathscr H }{\mathscr
H^2}-\frac{\mathscr W}{\mathscr H}\\\nonumber &=&\frac{\left(\p\bp
g_{i\bar j}\right) F^i\bar{F^j}}{\mathscr H}-\frac{g_{i\bar j}
T^i\wedge \bar
{T^j}}{\mathscr H}\\
&=&-\frac{ R_{i\bar j
k\bar\ell}f^i_{\alpha}\bar{f^j_{\beta}}f_\mu^k\bar{f_{\nu}^\ell}W^\mu\bar
W^\nu dz^\alpha\wedge d\bar z^\beta}{\mathscr H}\label{E4}
\end{eqnarray} where the last identity follows from the facts that \beq
\p\bp g_{i\bar j}=\frac{\p^2 g_{i\bar j}}{\p z^k\bar z^\ell}\p
f^k\wedge \bar{\p f^\ell},\ \ \ g_{i\bar j} T^i \bar{T^j}=g^{p\bar
q}\frac{\p g_{i\bar q}}{\p z^k}\frac{\p g_{p\bar j}}{\p\bar z^\ell}
F^i\bar {F^j}\p f^k\wedge \bar{\p f^\ell} \eeq and the curvature
formula of the Chern connection on $(N,g)$ \beq R_{ k\bar\ell i\bar
j}=-\frac{\p^2 g_{i\bar j}}{\p z^k\p\bar  z^\ell}+g^{p\bar
q}\frac{\p g_{i\bar q}}{\p z^k}\frac{\p g_{p\bar j}}{\p\bar
z^\ell}.\eeq

\noindent  Hence, by (\ref{E0}) and (\ref{E4}), we get $$ \sq\p\bp
\mathscr Y=\left(\sq \p\bp\log \mathscr H^{-1}\right)\cdot \mathscr
Y+\frac{\sq \mathscr W}{\mathscr H}-\frac{\sq R_{i\bar j
k\bar\ell}f^i_{\alpha}\bar{f^j_{\beta}}f_\mu^k\bar{f_{\nu}^\ell}W^\mu\bar
W^\nu dz^\alpha\wedge d\bar z^\beta}{\mathscr H}.$$ On the other
hand, by the definition equation (\ref{E2}) of $\mathscr W$, it is
easy to see that $\sq\mathscr W$ is a semi-positive $(1,1)$-form on
$\P(T_M)$. Hence, we obtain Theorem \ref{main0}.\qed

\vskip 1\baselineskip

For the conformal energy density (\ref{density02}), we also have a
similar inequality as in (\ref{S1}). \btheorem\label{variant} Let
$f:(M,h)\>(N,g)$ be a holomorphic map. Then for any $\phi\in
C^\infty(\P(T_M),\R)$, we have \beq \sq\p\bp \mathscr Y_\phi\geq
\left(\sq \p\bp\log {\mathscr H_\phi}^{-1}\right)\cdot {\mathscr
Y_\phi}-\frac{\sq R_{i\bar j
k\bar\ell}f^i_{\alpha}\bar{f^j_{\beta}}f_\mu^k\bar{f_{\nu}^\ell}W^\mu\bar
W^\nu dz^\alpha\wedge d\bar z^\beta}{ {\mathscr H_\phi}},\label{S03}
\eeq where $\mathscr Y_\phi =e^\phi \mathscr Y$ and  $\mathscr
H_\phi=\mathscr H e^{-\phi}$. \etheorem \bproof We use the same
notions as in the proof of Theorem \ref{main0}. It is easy to see
that $\mathscr F=\mathscr Y\cdot \mathscr H=\mathscr Y_\phi\cdot
\mathscr H_\phi$. By defining a new quantity as in the equation
(\ref{E2}) \beq \mathscr W_\phi=g_{i\bar j}\left(\p
F^i+F^i\p\log\mathscr H_\phi^{-1}+T^i\right)\wedge\bar{ \left(\p
F^j+F^j\p\log\mathscr H_\phi^{-1}+T^j\right)}, \eeq one can use
similar computations as in the proof of Theorem \ref{main0} to
deduce (\ref{S03}). \eproof

\vskip 1\baselineskip

\noindent\emph{The proof of  Theorem \ref{main}.}  We compute the
curvature form  $\sq\p\bp\log \mathscr H^{-1}$ over the projective
bundle $X:=\P(T_M)$ where $\pi:X\>M$ the natural projection,
following the  calculations in \cite[Proposition~9.2]{Gri65}.
Indeed, $\sq\p\bp\log \mathscr H^{-1}$ is the curvature of the
tautological line bundle $\left(\sO_{T^*_M}(-1),\mathscr H\right)$.
On the \emph{Hermitian manifold} $M$, we can choose ``normal
coordinates'' $\{z^\alpha\}_{\alpha=1}^m$ on a small open set
$U\subset M$ centered at point $p\in M$ such that \beq
h_{\alpha\bar\beta}(p)=\delta_{\alpha\bar\beta},\ \ \ \ \frac{\p
h_{\alpha\bar\beta}}{\p z^\gamma}(p)=-\frac{\p
h_{\gamma\bar\beta}}{\p z^\alpha}(p). \label{normal}\eeq By using
this local trivialization chart $U$ on $M$, we know $T_M|_U\cong
U\times \C^m$ and $\pi^{-1}(U)\cong U\times \P^{m-1}$. Let $q\in X$
such that $\pi(q)=p$. Locally, we write $p=(z_o^1,\cdots, z_o^m)$
and $q=(z_o^1,\cdots, z_o^m,[W_o^1,\cdots, W_o^m])$ where
$[W^1,\cdots, W^m]$ are the homogeneous coordinates on the fiber
$\P(T_pM)\cong \P^{m-1}$. Let $u=(W_o^1,\cdots, W_o^m,
\underbrace{0,\cdots, 0}_{m-1}) \in T_qX$. We claim that at point
$q\in X$, $u$ is a positive direction of $\sq\p\bp\log \mathscr
H^{-1}$, i.e. \beq \left(\sq\p\bp\log \mathscr H^{-1}\right)(u,\bar
u)>0. \eeq

\noindent  Indeed, if we set $H^{\gamma\bar\delta}=W^\gamma\bar
W^\delta$, then $\mathscr H=h_{\alpha\bar \beta}
H^{\alpha\bar\beta}$ and \be \p\bp\log\mathscr
H^{-1}&=&\frac{\p\mathscr H\wedge \bp\mathscr H}{\mathscr
H^2}-\frac{\p\bp \mathscr H}{\mathscr
H}\\
&=&\frac{\left(H^{\alpha\bar\delta}\p h_{\alpha\bar\beta}
+h_{\alpha\bar\beta}\p H^{\alpha\bar\beta}\right)\wedge
\left(H^{\gamma\bar\delta}\bp
h_{\gamma\bar\delta}+h_{\gamma\bar\delta}\p H^{\gamma\bar\delta}
\right)}{\mathscr H^2}\\
&&-\frac{\p\bp h_{\gamma\bar\delta}\cdot H^{\gamma\bar\delta}+\p
h_{\gamma\bar\delta}\wedge \bp H^{\gamma\bar\delta}+\p
H^{\gamma\bar\delta}\wedge \bp
h_{\gamma\bar\delta}+h_{\gamma\bar\delta}\p\bp
H^{\gamma\bar\delta}}{\mathscr H}.
 \ee

\noindent Note that $\p h_{\alpha\bar\beta}=\frac{\p
h_{\alpha\bar\beta}}{\p
 z^\mu}dz^\mu$ is along the base directions and $\p H^{\gamma\bar\delta}=\p \left(W^\gamma\bar W^\delta\right)$ is along the fiber
 directions.
 We evaluate the $(1,1)$ form $\sq\p\bp\log\mathscr
H^{-1}$ at point $q=(z_o^1,\cdots, z_o^m,[W_o^1,\cdots, W_o^m])$
with $u=(W_o^1,\cdots, W_o^m, 0,\cdots, 0) \in T_qX$. Since
$$\left(\p h_{\alpha\bar\beta} \wedge \bp H^{\gamma\bar\delta}\right)(u,\bar u)=0,\ \ \ (\p\bp H^{\alpha\bar\beta})(u,\bar u)=0,\ \ \ \left(\p H^{\alpha\bar\beta} \wedge \bp H^{\gamma\bar\delta}\right)(u,\bar u)=0,$$
we obtain \begin{eqnarray} \nonumber && \left(\p\bp\log
\mathscr H^{-1}\right)(u,\bar u)\\
 \nonumber&=&\left(-\frac{H_o^{\gamma\bar\delta}\p\bp
h_{\gamma\bar\delta}}{\mathscr
H_o}+\frac{H_o^{\alpha\bar\beta}H_o^{\gamma\bar\delta} \p
h_{\alpha\bar\beta}\wedge \bp h_{\gamma\bar\delta}}{\mathscr
H_o^2}\right)(u,\bar u)\\
\label{HSC}&=&\left(-\frac{\p^2 h_{\gamma\bar\delta}}{\p z^\mu\p\bar
z^\nu}+h^{\alpha\bar\beta}\frac{\p h_{\gamma\bar\beta}}{\p
z^\mu}\frac{\p h_{\alpha\bar\delta}}{\p \bar z^\nu}\right)
\frac{W_o^\mu\bar W_o^\nu W_o^\gamma\bar W_o^\delta }{|W_o|^2}\\
 \nonumber&=&\frac{R_{\mu\bar\nu \gamma\bar\delta} W_o^\mu\bar W_o^\nu
W_o^\gamma\bar W_o^\delta }{|W_o|^2}\end{eqnarray} where in the
second identity we used
the anti-symmetric property of the first order derivatives given in the normal coordinates (\ref{normal}).\\

  Suppose ${\mathscr Y}$ is not identically zero on $X$. Let $q\in
X$ be a maximum value point of $ {\mathscr Y}$ and $p=\pi(q)\in M$.
Hence $\mathscr Y(q)>0$ and $\sq\p\bp\mathscr Y(q)\leq 0$. By using
the previous setting, we can choose ``normal coordinates''
$\{z^\alpha\}_{\alpha=1}^m$ centered at point $p$ (e.g.
\ref{normal}) and we write $p=(z_o^1,\cdots, z_o^m)$ and
$q=(z_o^1,\cdots, z_o^m,[W_o^1,\cdots, W_o^m])$ where
$[W_o^1,\cdots, W^m_o]$ is the homogeneous coordinate on the fiber
$\P(T_pM)$.  Evaluating the inequality (\ref{S1})  on the vector
 $u=(W_o^1,\cdots, W_o^m, 0,\cdots, 0) \in T_qX$, we have
\beq 0\geq \left(\p\bp \mathscr  Y\right)(u,\bar u)\geq
\left(\p\bp\log \mathscr H^{-1}\right)(u,\bar u)\cdot \mathscr
Y-\frac{ R_{i\bar j
k\bar\ell}f^i_{\alpha}\bar{f^j_{\beta}}f_\mu^k\bar{f_{\nu}^\ell}W^\alpha_o\bar
W^\beta_oW_o^\mu\bar W_o^\nu }{\mathscr H}.\label{S4}\eeq

\noindent If we set $F^i=f^i_\alpha W^\alpha_o$, then by (\ref{S4})
and (\ref{HSC}), we obtain \beq 0\geq \frac{R_{\mu\bar\nu
\gamma\bar\delta} W_o^\mu\bar W_o^\nu W_o^\gamma\bar W_o^\delta
}{|W_o|^2}\cdot \mathscr Y(q)-\frac{ R_{i\bar j k\bar\ell} F^i\bar
F^j F^k\bar F^\ell}{\mathscr H}.\label{S5}\eeq Note that $\mathscr
Y(q)=\frac{g_{i\bar j} F^i\bar F^j}{\mathscr H}>0$ and so $(F^i)\neq
0$. If $(M,h)$ has positive (resp. nonnegative) holomorphic
sectional curvature and $(N,g)$ has non-positive (resp. negative)
holomorphic sectional curvature, then the right hand side of
(\ref{S5}) is positive, which is absurd. Therefore, ${\mathscr Y}$
must be identically zero on $X$ and so $f$ is a constant map.
 \qed

\vskip 1\baselineskip

\noindent\emph{The proof of Corollary \ref{main1}.} Fix an arbitrary
smooth metric $h$ on $M$ and let $\mathscr
H=h_{\gamma\bar\delta}W^\gamma\bar W^\delta$. Suppose
$\sO_{T_M^*}(-1)$ is RC-positive, then there exists a Hermitian
metric $\tilde {\mathscr H}$ on $\sO_{T_M^*}(-1)$ such that its
curvature $-\sq\p\bp\log \tilde {\mathscr H}$ has at least one
positive eigenvalue at each point of $\P(T_M)$. Since $\mathscr H$
is also a smooth Hermitian metric on $\sO_{T_M^*}(-1)$, there exists
$\phi\in C^\infty(\P(T_M),\R)$ such that ${\tilde {\mathscr
H}}=\mathscr H e^{-\phi}$ and \beq -\sq\p\bp\log \tilde {\mathscr
H}=-\sq\p\bp\log\mathscr H+\sq\p\bp\phi. \eeq  We consider $\tilde
{\mathscr Y}=\mathscr Ye^{\phi}$ on $\P(T_M)$.  By Theorem
\ref{variant}, we obtain \beq \sq\p\bp \tilde{\mathscr Y}\geq
\left(\sq \p\bp\log \tilde {\mathscr H}^{-1}\right)\cdot
\tilde{\mathscr Y}-\frac{\sq R_{i\bar j
k\bar\ell}f^i_{\alpha}\bar{f^j_{\beta}}f_\mu^k\bar{f_{\nu}^\ell}W^\mu\bar
W^\nu dz^\alpha\wedge d\bar z^\beta}{\tilde {\mathscr H}}. \eeq
 Suppose $\tilde
{\mathscr Y}$ is not identically zero on $\P(T_M)$. Let $p\in
\P(T_M)$ be a maximum value point of $\tilde {\mathscr Y}$. Hence,
$\tilde{\mathscr Y}(p)>0$ and  \beq \sq\p\bp \tilde {\mathscr
Y}(p)\leq 0.\label{E5} \eeq  On the other hand, since $(N,g)$ has
non-positive holomorphic bisectional curvature, \beq  \sq\p\bp
\tilde{\mathscr Y}\geq \left(\sq \p\bp\log \tilde {\mathscr
H}^{-1}\right)\cdot \tilde{\mathscr Y}. \eeq This contradicts to
(\ref{E5}) since $\sq \p\bp\log \tilde {\mathscr H}^{-1}$ has at
least one positive eigenvalue at point $p$ and $\tilde{\mathscr
Y}(p)>0$.\qed

\vskip 1\baselineskip

\noindent\emph{The sketched proof of Corollary \ref{main2}.} Fix an
arbitrary smooth metric $h$ on $M$ and set $\mathscr
H=h_{\gamma\bar\delta}W^\gamma\bar W^\delta$. Suppose
$\sO_{T_M^*}(-1)$ is RC-positive, then there exists   $\phi\in
C^\infty(\P(T_M),\R)$ such that $\mathscr H e^{-\phi}$ is a smooth
Hermitian metric on $\sO_{T_M^*}(-1)$ and \beq -\sq\p\bp\log\mathscr
H+\sq\p\bp\phi \eeq has at least one positive eigenvalue at each
point in $\P(T_M)$.

Since $\sO_{T^*N}(1)$ is nef, by \cite[Theorem~1.12]{DPS94}, for any
given Hermitian metric $\omega_N$ on $N$, there exist Hermitian
metrics $g_k$ on $\mathrm{Sym}^{\ts k}T_N$ and $\eps_k>0$ such that
\beq R_{i\bar j K\bar L} u^i\bar u^j V^K\bar V^\ell\leq k\eps_k
|u|_{\omega_N}^2|V|_{g_k}^2\label{nefe}\eeq where $\eps_k\>0$ as
$k\>\infty$.
 Let $\{W^A\}$ be the  symmetric polynomials in
$\{W^\alpha\}_{\alpha=1}^m$ of degree $k$. We can define the higher
order energy density on $\P(T_M)$ as following:

\beq  \mathscr Y_k=g_{I\bar J}F_A^I\bar{F}_{B}^J\frac{W^A\bar
W^B}{\left(\sum h_{\gamma\bar\delta}W^\gamma \bar
W^{\delta}\right)^k}, \label{density04}\eeq where the matrix
$\left(F^I_A\right)$ is the $k$-th symmetric power product
$\mathrm{Sym}^{\ts k}\left(f^i_\alpha\right)$ and $\left(g_{I\bar
J}\right)$ is the metric matrix of $g_k$ on $\mathrm{Sym}^{\ts
k}T_N$ with respect to the given trivialization on $N$. We use the
conformal change (\ref{density02}) of the generalized energy density
$\mathscr Y_k$, i.e. \beq \mathscr Y_{k,\phi}=e^{k\phi}\mathscr Y_k.
\eeq
 As in Theorem  \ref{main0} and
Theorem \ref{variant}, it is not hard to deduce the inequality \beq
\sq\p\bp \mathscr Y_{k,\phi}\geq \left(\sq \p\bp\log \mathscr
H_\phi^{-k}\right)\cdot \mathscr Y_{k,\phi}-\frac{\sq R_{i\bar j
K\bar L}f^i_{\alpha}\bar{f^j_{\beta}}F_C^K\bar{F_{D}^L}W^C\bar W^D
dz^\alpha\wedge d\bar z^\beta}{\mathscr H_\phi^k},\label{E} \eeq
where $R_{i\bar j K\bar L}$ is the curvature tensor component of the
Hermitian vector bundle $\left(\mathrm{Sym}^{\ts k}T_N,g_k\right)$.
 Since $M$ is compact and $\left(\sO_{T^*_M}(-1),\mathscr
{H}e^{-\phi}\right)$ is RC-positive, we deduce \beq \min_{P\in
\P(T_M)}\sup_{u\in T_P \P(T_M), u\neq 0}\frac{\left(\sq \p\bp\log
\mathscr H_\phi^{-1}\right)(u, \bar u)}{|u|_\omega^2}\geq 4C
\label{RCestimate}\eeq where $\omega$ is a fixed metric on $\P(T_M)$
and $C=C(\omega, \mathscr H_{\phi})$ is a positive constant.

 Suppose ${\mathscr Y_{k,\phi}}$ is
not identically zero on $\P(T_M)$. Let $P\in \P(T_M)$ be a maximum
value point of ${\mathscr Y_{k,\phi}}$. Hence,   $ \sq\p\bp
{\mathscr Y_{k,\phi}}(P)\leq 0$. Let $u\in T_P(\P(T_M))$ be a
positive direction of $\sq \p\bp\log \mathscr H_\phi^{-1}$ such that
\beq \left(\sq \p\bp\log \mathscr H_\phi^{-1}\right)(u,\bar u)\geq
2C|u|^2_{\omega}.\label{RCest}\eeq

\noindent By (\ref{E}), (\ref{nefe}) and (\ref{RCest}), we deduce
\beq \left( \sq\p\bp \mathscr Y_{k,\phi}\right)(u,\bar u)\geq
2kC|u|^2_{\omega}\mathscr Y_{k,\phi}- k\eps_k \cdot |\p
f|^2_{\omega\ts f^*\omega_N}\cdot |u|^2_{\omega}\cdot\mathscr
Y_{k,\phi}.\eeq Since $\eps_k\>0$ as $k\>\infty$, we deduce that
when $k$ is large enough \beq 2C-\eps_k \cdot |\p f|^2_{\omega\ts
f^*\omega_N}>C.\eeq Hence, we obtain $$\left( \sq\p\bp \mathscr
Y_{k,\phi}\right)(u,\bar u)\geq kC|u|^2_{\omega}\mathscr
Y_{k,\phi}>0$$ which is absurd.  We conclude $\mathscr Y_{k,\phi}$
is identically zero and so $f$ is a constant map.\qed

\vskip 1\baselineskip

Corollary \ref{main3} can be proved by using similar strategies as
in the proof of Corollary \ref{main2}, and we need to use the
generalized version of the energy density (\ref{density2}) on
$\P(f^*T_N^*)$. The negative holomorphic sectional curvature is used
in a similar way as in the proof Theorem \ref{main}.

By using the proof in  Corollary \ref{main2}, we get the following
generalization.

\bcorollary\label{main22} Let $f:M\>N$ be a holomorphic map between
compact complex
 manifolds. If $\sO_{\mathrm{Sym}^{\ts k}T^*_M}(-1)$ is RC-positive for some $k\geq 1$ and $\sO_{T^*_N}(1)$ is nef, then $f$ is a constant map.
\ecorollary

\vskip 2\baselineskip

\section{The Hessian estimates and rigidity of harmonic maps and pluri-harmonic maps}\label{4}

 In this section,  we use the Hessian estimate of the generalized energy density to investigate harmonic maps and pluri-harmonic from complex manifolds to Riemannian manifolds and K\"ahler manifolds, and establish Theorem \ref{main10}, Theorem \ref{main11},
 Theorem \ref{main12} and Theorem \ref{main13}.

\subsection{The Hessian estimate and rigidity of pluri-harmonic maps into
Riemannian manifolds}

\noindent Now we are ready to prove Theorem \ref{main10}, i.e. a
projective bundle version of Lemma \ref{cri}.\\

\noindent\emph{The proof of Theorem \ref{main10}.} We shall use
similar strategies as in the proof of Theorem \ref{main0}. We choose
an arbitrary point $p\in M$ and $q=f(p)\in N$. Let $\{z^\al\}$ and
$\{x^i\}$ be coordinates centered at point $p\in M$ and $q\in N$
respectively such that
\begin{eqnarray*}h_{\alpha\bar\beta}(p)=\delta_{\alpha\bar\beta},\ \ \  g_{ij}(q)=\delta_{ij},\ \ \  \frac{\p g_{ij}}{\p x^\ell}(q)=0.\end{eqnarray*}
Since $f$ is pluri-harmonic, by formula (\ref{pluriharmonic}),  at
point $p$, we have \beq f^i_{\alpha\bar\beta}(p)=\frac{\p^2 f^i}{\p
z^\alpha\bar z^\beta}(p)=0. \label{plurinormal}\eeq

\noindent We write $\mathscr F=\mathscr Y\cdot \mathscr H$. A
straightforward calculation on $X:=\P(T_M)$ yields \beq
\p\bp\mathscr Y=\left(\p\bp\log \mathscr H^{-1}\right) \mathscr
Y+\frac{\p\bp \mathscr F}{\mathscr H}+ \frac{\mathscr F\p \mathscr
H\wedge \bp \mathscr H}{\mathscr H^3}-\frac{\p \mathscr H\wedge \bp
\mathscr F+\p\mathscr F\wedge \bp \mathscr H }{\mathscr H^2}.  \eeq
We write
 $F^i=f^i_\alpha
 W^\alpha$ and  $\mathscr F=g_{i j}F^i\bar F^j$.  Hence,
$$\p \mathscr F=\p g_{ij} \cdot F^i\bar{F^j}+g_{i j}\cdot \p
F^i\cdot F^j+g_{ij}F^i \p f^i_{\bar\beta} \cdot \bar W^\beta,$$ and
$$ \bp\mathscr F=\bp g_{i j}\cdot F^i\bar{F^j}+ g_{i j} F^i\bar{\p
F^j}+g_{ij}\bar F^j \bp f^i_\alpha \cdot W^\alpha. $$ By using
(\ref{plurinormal}), at point $P\in X$ with $\pi(P)=p$, we have
$$\p \mathscr F=\p g_{ij} \cdot F^i\bar{F^j}+g_{i j}\cdot \p
F^i\cdot F^j,\ \ \  \bp\mathscr F=\bp g_{i j}\cdot F^i\bar{F^j}+
g_{i j} F^i\bar{\p F^j}. $$

\noindent  By a similar computation as in (\ref{E3}), we deduce
\begin{eqnarray} \nonumber&&\frac{\p\bp \mathscr F}{\mathscr H}+ \frac{\mathscr F\p
\mathscr H\wedge \bp \mathscr H}{\mathscr H^3}-\frac{\p \mathscr
H\wedge \bp \mathscr F+\p\mathscr F\wedge \bp \mathscr H }{\mathscr
H^2}\\ \nonumber&=&\frac{\p\bp g_{i j}\cdot F^i\bar{F^j}+ \p g_{i
j}\cdot F^i\cdot \bar{\p F^j}-\bp g_{i j}\cdot\p
F^i\cdot \bar{F^j}+ g_{i j}\p F^i\bar{\p F^j}}{\mathscr H}\\ \nonumber&&+\frac{g_{ij}\cdot \p\bp F^i \cdot \bar F^j+g_{ij}\cdot F^i\cdot\p\bp\bar F^j}{\mathscr H^2}\\
\nonumber&&+\frac{\mathscr F \p\log\mathscr H^{-1}\wedge
\bp\log\mathscr H^{-1}}{\mathscr H}\\
\nonumber&&-\frac{\p\log\mathscr H^{}\cdot \left(\bp g_{i j}\cdot
F^i\bar{F^j}+ g_{i j} F^i\bar{\p F^j}\right)+\left(\p g_{i j} \cdot
F^i\bar{F^j}+g_{i j}\cdot \p F^i\cdot F^j\right)\cdot
\bp\log\mathscr H^{}}{\mathscr H}.\end{eqnarray}

 \noindent The key ingredient
is to define the $(1,1)$-form $\mathscr W$ on $X$ as in (\ref{E2}):
\beq \mathscr W:=g_{i j}\left(\p F^i+F^i\p\log\mathscr
H^{-1}+T^i\right)\wedge\bar{ \left(\p F^j+F^j\p\log\mathscr
H^{-1}+T^j\right)} \label{E8} \eeq where $T^i=\Gamma_{pk}^iF^k\p
f^p$. By using a similar computation as in (\ref{E4}), we obtain
\begin{eqnarray}&&\nonumber\frac{\p\bp \mathscr F}{\mathscr H}+ \frac{\mathscr F\p
\mathscr H\wedge \bp \mathscr H}{\mathscr H^3}-\frac{\p \mathscr
H\wedge \bp \mathscr F+\p\mathscr F\wedge \bp \mathscr H }{\mathscr
H^2}-\frac{\mathscr W}{\mathscr H}\\\nonumber &=&\frac{\left(\p\bp
g_{i j}\right) F^i\bar{F^j}}{\mathscr H}-\frac{g_{i j} T^i\wedge
\bar
{T^j}}{\mathscr H}+\frac{g_{ij}\cdot \p\bp F^i \cdot \bar F^j+g_{ij}\cdot F^i\cdot\p\bp\bar F^j}{\mathscr H^2}\\
\nonumber &=&\frac{\left(\p\bp g_{i j}\right) F^i\bar{F^j}}{\mathscr
H}-\frac{g_{i j} T^i\wedge \bar
{T^j}}{\mathscr H}+\frac{g_{ij}\cdot \left(\p\bp f_\gamma^i\right)W^\gamma \cdot \bar F^j
+g_{ij}\cdot F^i\cdot\left(\p\bp f^j_{\bar\delta}\right)\bar W^\delta}{\mathscr H^2}\\
\nonumber &=&\frac{\displaystyle{\left(\frac{\p^2 g_{ij}}{\p x^k\p
x^\ell}f^k_{\al}f^\ell_{\bbe}
f^i_{\ga}f^j_{\bar\delta}+f^i_{\al\bbe\ga}f^i_{\bar\delta}+f^i_{\ga}f^i_{\al\bbe\bar\delta}\right)}W^\gamma\bar
W^\delta dz^\alpha\wedge d\bar z^\beta}{\mathscr H}
\\\nonumber&=&-\frac{ \left(R_{i\ell k j}+R_{ik\ell
j}\right)f^i_{\alpha}{f^j_{\bar\beta}}f_\gamma^k{f_{\bar\delta}^\ell}W^\gamma\bar
W^\delta dz^\alpha\wedge d\bar z^\beta}{\mathscr H},
\end{eqnarray}
where the last identity follows from a standard computation  by
using the pluri-harmonic equation (\ref{pluriharmonic}). Indeed, by
(\ref{pluriharmonic}), we obtain
$$f^i_{\al\bbe\ga}f^i_{\bar\delta}=-\frac{\p\Ga^\ell_{ij}}{\p
x^k}f^i_{\al}f^j_{\bbe}f^k_{\ga}f^\ell_{\bar\delta},\ \ \
 f^i_{\al\bbe\bar\delta}f^i_{\ga}=-\frac{\p\Ga^k_{ij}}{\p x^\ell}f^i_{\al}f^j_{\bbe}f^k_{\ga}f^\ell_{\bar\delta},$$
and
$$\frac{\p^2 g_{ij}}{\p x^k\p
x^\ell}f^k_{\al}f^l_{\bbe}
f^i_{\ga}f^j_{\bar\delta}+f^i_{\al\bbe\ga}f^i_{\bar\delta}+f^i_{\ga}f^i_{\al\bbe\bar\delta}=\left(\frac{\p^2
g_{k\ell}}{\p x^i\p x^j}-\frac{\p\Ga^\ell_{ij}}{\p
x^k}-\frac{\p\Ga^k_{ij}}{\p
x^\ell}\right)f^i_{\al}f^j_{\bbe}f^k_{\ga}f^\ell_{\bar\delta}.$$

\noindent By using the Riemannian curvature tensor
(\ref{rcurvature}), it is easy to show that
 \beq \frac{\p^2 g_{k\ell}}{\p x^i\p
x^j}-\frac{\p\Ga^\ell_{ij}}{\p x^k}-\frac{\p\Ga^k_{ij}}{\p x^\ell}=-
\left(R_{i\ell k j}+R_{ik\ell j}\right).\label{key3}\eeq
 By using the constraint equation (\ref{D}),
we get \be && \frac{\p\bp \mathscr F}{\mathscr H}+ \frac{\mathscr
F\p \mathscr H\wedge \bp \mathscr H}{\mathscr H^3}-\frac{\p \mathscr
H\wedge \bp \mathscr F+\p\mathscr F\wedge \bp \mathscr H }{\mathscr
H^2}-\frac{\mathscr W}{\mathscr H}\\&=&-\frac{ R_{i\ell k j}
f^i_{\alpha}{f^j_{\bar\beta}}f_\gamma^k{f_{\bar\delta}^\ell}W^\gamma\bar
W^\delta dz^\alpha\wedge d\bar z^\beta}{\mathscr H}.\ee Finally, we
obtain
$$ \sq\p\bp
\mathscr Y=\left(\sq \p\bp\log \mathscr H^{-1}\right)\cdot \mathscr
Y+\frac{\sq \mathscr W}{\mathscr H}-\frac{\sq R_{i\ell k j}
f^i_{\alpha}{f^j_{\bar\beta}}f_\gamma^k{f_{\bar\delta}^\ell}W^\gamma\bar
W^\delta dz^\alpha\wedge d\bar z^\beta}{\mathscr H}$$  and the proof
of Theorem \ref{main10} is completed.\qed

\vskip 1\baselineskip

\noindent\emph{The proof of Theorem \ref{main11}.} The proof is
similar to that in Corollary \ref{main1}.  Fix an arbitrary smooth
metric $h$ on $M$ and let $\mathscr
H=h_{\gamma\bar\delta}W^\gamma\bar W^\delta$.  Since
$\sO_{T_M^*}(-1)$ is RC-positive,  there exist a Hermitian metric
$\tilde {\mathscr H}$ on $\sO_{T_M^*}(-1)$ and $\phi\in
C^\infty(\P(T_M),\R)$ such that ${\tilde {\mathscr H}}=\mathscr H
e^{-\phi}$ and \beq -\sq\p\bp\log \tilde {\mathscr
H}=-\sq\p\bp\log\mathscr H+\sq\p\bp\phi \eeq has at least one
positive eigenvalue at each point in $\P(T_M)$. We consider $\tilde
{\mathscr Y}=\mathscr Ye^{\phi}$ on $\P(T_M)$. By using similar
computations as in the proof of Theorem \ref{main10} and Theorem
\ref{variant}, we obtain \beq  \sq\p\bp \tilde{\mathscr Y}\geq
\left(\sq \p\bp\log \tilde {\mathscr H}^{-1}\right)\cdot
\tilde{\mathscr Y}-\frac{\sq R_{i\ell k j}
f^i_{\alpha}{f^j_{\bar\beta}}f_\gamma^k{f_{\bar\delta}^\ell}W^\gamma\bar
W^\delta dz^\alpha\wedge d\bar z^\beta}{\tilde {\mathscr H}}.
\label{E9} \eeq

\noindent
 Suppose $\tilde
{\mathscr Y}$ is not identically zero on $\P(T_M)$. Let $p\in
\P(T_M)$ be a maximum value point of $\tilde {\mathscr Y}$. Hence,
$\tilde{\mathscr Y}(p)>0$ and  \beq \sq\p\bp \tilde {\mathscr
Y}(p)\leq 0.\label{E10} \eeq  Since $\sq \p\bp\log \tilde {\mathscr
H}^{-1}$ has at least one positive eigenvalue at point $p$, there
exists a non-zero  vector $u=(a_1,\cdots, a_m, v_1,\cdots, v_{m-1})$
such that \beq \sq \p\bp\log \tilde {\mathscr H}^{-1}(u,\bar u)>0.
\label{E11}\eeq 

\noindent Let $H^{\gamma\bar\delta}=W^\gamma\bar W^\delta$ and  \beq
C_{\alpha\bar\beta}=R_{i\ell k j}f^i_{\al}
f^j_{\bar\beta}\left(H^{\gamma\bar\delta}f^k_{\ga}
f^\ell_{\bar\delta}\right).
 \eeq
By using the constraint equation (\ref{D})
  for pluri-harmonic maps, we have
   \beq C_{\alpha\bar\beta}=R_{i\ell k j}f^i_{\al}f^j_{\bbe}
\left(H^{\gamma\bar\delta}f^k_{\ga}f^\ell_{\bar\delta}+H^{\gamma\bar\delta}f^\ell_{\ga}f^k_{\bar\delta}\right)
\eeq On the other hand, since  $\sq
C_{\alpha\bar\beta}dz^\alpha\wedge d\bar z^\beta$ is a real $(1,1)$
form, we obtain \beq
C_{\alpha\bar\beta}=\bar{C_{\beta\bar\alpha}}=R_{i\ell k
j}f^j_{\al}f^i_{\bbe}
\left(H^{\gamma\bar\delta}f^k_{\ga}f^\ell_{\bar\delta}+H^{\gamma\bar\delta}f^\ell_{\ga}f^k_{\bar\delta}\right).
 \eeq
Therefore,
 \beq C_{\alpha\bar\beta}=\frac{1}{2}R_{i\ell k
j}\left(f^i_{\al}f^j_{\bbe}+f^j_{\al}f^i_{\bbe}\right)
\left(H^{\gamma\bar\delta}f^k_{\ga}f^\ell_{\bar\delta}+H^{\gamma\bar\delta}f^\ell_{\ga}f^k_{\bar\delta}\right).
\eeq  If we set $u^i=\sum_\alpha f^i_\alpha a_\alpha$ and
$F^k=\sum_\gamma f^k_\gamma W^\gamma$, then at point $p$, \beq
\sum_{\alpha,\beta} C_{\alpha\bar\beta}a_\alpha\bar a_\beta
=\frac{1}{2}\sum R_{i\ell k j}\left(u^i \bar u^j+u^j\bar
u^i\right)\left(F^k \bar F^\ell+F^\ell \bar F^k\right).\eeq Let
$u^i=c^i+\sq b^i$ and $F^i=A^i+\sq B^i$ where $c^i, b^i, A^i, B^i$
are real numbers. Therefore we have
 \beq
\sum_{\alpha,\beta} C_{\alpha\bar\beta}a_\alpha\bar a_\beta =2\sum
R_{i\ell k j}\left(c^i c^j+b^i b^j\right)(A^k A^\ell+B^k
B^\ell).\eeq Since $(N,h)$ has non-positive Riemannian sectional
curvature, we deduce \beq \sum_{\alpha,\beta}
C_{\alpha\bar\beta}a_\alpha\bar a_\beta =2\sum_{}R_{i\ell k
j}\left(c^i c^j+b^i b^j\right)(A^k A^\ell+B^k B^\ell)\leq 0.
\label{nonnegative}\eeq

\noindent
 By formulas (\ref{E9}),
 (\ref{E10}), (\ref{E11}) and (\ref{nonnegative}), we get a contradiction.
 Hence $\tilde{\mathscr Y}$ is a constant and $f$ must be a constant map. \qed

\vskip 1\baselineskip

\noindent By using  formula (\ref{hessian2}) and the conformal
change technique, we also obtain:

\btheorem \label{main4} Let $f:M\>(N,g)$ be a pluri-harmonic map
from a  compact complex manifold to a Riemannian manifold $(N,g)$
with non-positive Riemannian sectional curvature. If there exist a
Hermitian metric $\omega$ on $M$ and a Hermitian metric $h$ on
$T^{1,0}M$ such that
$$\mathrm{tr}_\omega R^{(T^{1,0}M,h)}\in\Gamma(M,\mathrm{End}(T^{1,0}M))$$ is quasi-positive, then $f$ is
a constant map. \etheorem

\bproof Let $\omega_G=e^f\omega$ be a smooth Gauduchon metric in the
conformal class of $\omega$, i.e. $\p\bp\omega_G^{m-1}=0$. Let
$\omega_{G}=\sq G_{\alpha\bar\beta}dz^\alpha\wedge d\bar z^\beta$.
By taking trace of (\ref{hessian}), we obtain
 \begin{eqnarray} \mathrm{tr}_{\omega_G}\sq\p\bp u\nonumber\geq
\left(G^{\alpha\bar\beta}R^h_{\alpha\bar \beta \gamma\bar
\delta}\right) h^{\gamma \bar \nu} h^{\mu\bar \delta} g_{ij} f^{i
}_\mu {f^{j}_{\bar{\nu}}}  - R_{i\ell k
j}\left(G^{\alpha\bar\beta}f^i_{\al}f^j_{\bbe}\right)
\left(h^{\gamma\bar\delta}f^k_{\ga}f^\ell_{\bar\delta}\right).
\label{conformal}\end{eqnarray} By using a similar argument as
above,  we have
$$R_{i\ell k j}\left(G^{\alpha\bar\beta}f^i_{\al}f^j_{\bbe}\right)
\left(h^{\gamma\bar\delta}f^k_{\ga}f^\ell_{\bar\delta}\right)=
\frac{1}{2}R_{i\ell k
j}\left(G^{\alpha\bar\beta}f^i_{\al}f^j_{\bbe}+G^{\alpha\bar\beta}f^j_{\al}f^i_{\bbe}\right)
\left(h^{\gamma\bar\delta}f^k_{\ga}f^\ell_{\bar\delta}+h^{\gamma\bar\delta}f^\ell_{\ga}f^k_{\bar\delta}\right).
$$ At a point $p\in M$, we can assume
$h_{\alpha\bar\beta}(p)=\delta_{\alpha\beta}$ and
$G^{\alpha\bar\beta}(p)=\lambda_\alpha\delta_{\alpha\beta}$ where
$\lambda_\alpha>0$. Hence,
$$R_{i\ell k j}\left(G^{\alpha\bar\beta}f^i_{\al}f^j_{\bbe}\right)
\left(h^{\gamma\bar\delta}f^k_{\ga}f^\ell_{\bar\delta}\right)=
2\sum_{\alpha,\gamma}R_{i\ell k j}\lambda_\alpha(A_\alpha^i
A_\alpha^j+B_\alpha^i B_\alpha^j)(A_\gamma^k
A_\gamma^\ell+B_\gamma^k B_\gamma^\ell)\leq 0$$ where
$f_\alpha^i=A_\alpha^i+\sq B_\alpha^i$ and $A_\alpha^i, B_\alpha^i$
are real numbers.

 On the other hand, since $\mathrm{tr}_\omega R^{(T^{1,0}M,h)}$ is
 quasi-positive, we know $$\mathrm{tr}_{\omega_G} R^{(T^{1,0}M,h)}=e^{-f}\mathrm{tr}_\omega R^{(T^{1,0}M,h)}$$
is also quasi-positive. Therefore,
$$\left(G^{\alpha\bar\beta}R^h_{\alpha\bar \beta \gamma\bar
\delta}\right) h^{\gamma \bar \nu} h^{\mu\bar \delta} g_{ij} f^{i
}_\mu {f^{j}_{\bar{\nu}}}\geq 0.$$ Since $\omega_G$ is Gauduchon, by
(\ref{conformal}) we deduce
$$\int_M \left(G^{\alpha\bar\beta}R^h_{\alpha\bar \beta \gamma\bar
\delta}\right) h^{\gamma \bar \nu} h^{\mu\bar \delta} g_{ij} f^{i
}_\mu {f^{j}_{\bar{\nu}}} \cdot\omega_G^m=0$$ and $\p f$ must be
identically zero on the open set where the curvature
$\left(G^{\alpha\bar\beta}R^h_{\alpha\bar \beta \gamma\bar
\delta}\right)$ is strictly positive. Since pluri-harmonic maps are
also Hermitian harmonic, by \cite[Theorem~6]{JY93}, $f$ must be a
constant map. \eproof

\subsection{Rigidity  of  pluri-harmonic maps into K\"ahler manifolds}

In this subsection, we shall prove  Theorem \ref{main12}. Let $M$ be
a complex manifold, $(N,g)$ be a Hermitian manifold and $f:M\>N$ be
a smooth map.  We denote by $\mathcal E=f^*(T^{1,0}N)$ and endow it
with the induced \textbf{Chern connection} $\nabla^\mathcal E$ from
$T^{1,0}N$. There is a natural decomposition $\nabla^{\mathcal
E}=\bp_\mathcal E+\p_\mathcal E$. Let $\{z^\alpha\}_{\alpha=1}^m$ be
the local holomorphic coordinates on $M$ and $\{w^i\}_{i=1}^n$ be
the local holomorphic coordinates on $N$. Let $e_i=f^*(\frac{\p}{\p
w^i})$. There are three $\mathscr E$-valued $1$-forms, i.e., \beq
\bp f =\frac{\p f^i}{\p \bar z^\beta}d\bar z^\beta\ts e_i, \ \ \p
f=\frac{\p f^i}{\p z^\alpha}dz^\alpha\ts e_i,\ \ \  df=\bp f+\p f.
\eeq

\bdefinition  A smooth map $f: M\>(N,g)$ from a complex manifold $M$
to a Hermitian manifold $(N,g)$  is called \emph{pluri-harmonic} if
it satisfies $\p_\mathcal E\bp f=0$, i.e. \beq \left(\frac{\p^2
f^i}{\p z^\alpha\p\bar z^\beta}+\Gamma_{jk}^i\frac{\p f^j}{\p
z^\alpha}\frac{\p f^k}{\p \bar z^\beta}\right)dz^\alpha\wedge d\bar
z^\beta\ts e_i=0,\label{pluridefinitioncc}\eeq where
$\Gamma_{jk}^i=g^{i\bar\ell}\frac{\p g_{k\bar\ell}}{\p z^j}$ is the
Christoffel symbol of the Chern connection on $(T^{1,0}N,g)$.

 \edefinition

\bremark If $(N,h)$ is not K\"ahler, in general, the indices
 $j$ and $k$  in the formula
 (\ref{pluridefinitioncc}) are not symmetric. In particular, we
do not have a constraint equation as in (\ref{D}) since neither
$\p_{\mathcal E} \p f=0$ nor $\bp_{\mathcal E}\bp f=0$ holds. Note
that the pluri-harmonic maps in (\ref{pluriharmonic}) and
(\ref{pluridefinitioncc}) are different since the metric connections
on the target are different.
 \eremark

\blemma\label{same} Let $f:M\>(N,g)$ be a pluri-harmonic map from a
complex manifold $M$ to a K\"ahler manifold $(N,g)$. Suppose
$g_{\R}$ is the background Riemannian metric of the K\"ahler metric
$g$ on $N$. Then $f:M\>(N,g_{\R})$ is a pluri-harmonic map in the
sense of (\ref{pluriharmonic}). \elemma \bproof It follows from the
standard complexification since the Chern connection
 is the same as the Levi-Civita connection  when the  manifold $(N,g)$ is K\"ahler. \eproof

\bremark\label{siure} If $f:M\>(N,g)$ is a pluri-harmonic map from a
complex manifold $M$ to a K\"ahler manifold $(N,g)$, we can get a
similar constraint equation as in (\ref{D}). Indeed, it is exactly
the complexification of (\ref{D}). Moreover, the complexification of
 the constraint term (\ref{hatC}) is
\beq  \sum_{\alpha,\gamma}R_{i\bar j k\bar \ell}\left(\frac{\p
f^i}{\p z^\alpha}\frac{\p \bar f^j}{\p z^\gamma}-\frac{\p f^i}{\p
z^\gamma}\frac{\p \bar f^j}{\p z^\alpha}\right)\bar{\left(\frac{\p
f^\ell}{\p z^\alpha}\frac{\p \bar f^k}{\p z^\gamma}-\frac{\p
f^\ell}{\p z^\gamma}\frac{\p \bar f^k}{\p z^\alpha}\right)}=0,\eeq
which is the notion introduced by Siu (\cite{Siu80}).  If $(N,g)$
has strongly negative curvature in the sense of Siu and
$\mathrm{rank}_{\R} df\geq 4$, then the pluri-harmonic map $f$ is
holomorphic or anti-holomorphic.
 \eremark

 \noindent \emph{The proof of Theorem \ref{main12}.} It follows from Lemma
\ref{same} and Theorem \ref{main11}. \qed

\vskip 1\baselineskip

\noindent \emph{The proof of Corollary \ref{coro}.} By using the
formula (\ref{HSC}), we deduce that if a Hermitian manifold $(M,h)$
has positive holomorphic sectional curvature, then
$(\sO_{T^*M}(-1),\mathscr H)$ is RC-positive. Hence, Corollary
\ref{coro} follows from Theorem \ref{main11}.\qed

\subsection{Rigidity  of  harmonic maps into Riemannian manifolds} In this subsection, we shall prove  Theorem \ref{main13}. Let $(M,h)$ be a compact Hermitian manifold, $(N,g)$ a Riemannian
manifold and $\mathscr E=f^*(TN)$ with the induced Levi-Civita
connection.
 $f$ is called \emph{Hermitian harmonic} if it satisfies $\mathrm{tr}_{\omega_h}\p_\mathscr E\bp
f=0$, i.e. \beq h^{\alpha\bar\beta}\left(\frac{\p^2 f^i}{\p
z^\alpha\p\bar z^\beta}+\Gamma_{jk}^i\frac{\p f^j}{\p \bar
z^\beta}\frac{\p f^k}{\p z^\alpha}\right)\ts e_i=0.\label{pseudo2}
\eeq \noindent It is easy to see that \bcorollary Let
$f:(M,h)\>(N,g)$ be a smooth map from a compact K\"ahler manifold
$(M,h)$ to a Riemannian manifold $(N,g)$. Then Hermitian harmonic
maps and harmonic maps coincide. \ecorollary

As analogous to Siu's strong negativity \cite{Siu80}, Sampson
 proposed  in \cite{Sam85} the following definition (see also
\cite[Definition~4.2]{JY91}):

\bdefinition \label{samdef}Let $(N,g)$ be a  Riemannian manifold.
The (complexified) curvature tensor $R$ of $(N,g)$ is said to have
\emph{non-positive complex sectional curvature} if \beq R(Z,\bar W,
W,\bar Z)\leq 0 \label{cscdef}\eeq for any $Z,W\in T_\C
N$.\edefinition

\noindent If $(N,g)$ has non-positive complex sectional curvature,
then it has non-positive Riemannian sectional curvature. Moreover,
 a K\"ahler manifold $(N,g)$ has
strongly non-positive curvature in the sense of Siu if and only if
its background Riemannian metric has non-positive complex sectional
curvature (e.g. \cite[Theorem~4.4]{LSYY17}).\\

By using Siu's $\p\bp$-trick and ideas in \cite{Siu80, Sam85, Sam86,
JY91, JY93}, one has the following result (see
\cite[Theorem~6.11]{LY14} and \cite[Theorem~4.2]{WY18}).

\blemma\label{pluri} Let $(M,h)$ be a compact  astheno-K\"ahler
manifold (i.e. $\p\bp\omega^{m-2}_h=0$) and $(N,g)$ a  Riemannian
manifold. Let $f:(M,h)\>(N,g)$ be a Hermitian harmonic map. If
$(N,g)$ has
 non-positive complex sectional  curvature, then $f$ is pluri-harmonic. \bproof If $f$
is Hermitian harmonic, i.e., $\mathrm{tr}_{\omega_h}\p_\mathscr E\bp
f=0$, it is easy to see that\beq \p\bp\{\bp f, \bp f\}
\frac{\omega_h^{m-2}}{(m-2)!}=4 |\p_\mathscr E\bp
f|^2\frac{\omega^m_h}{m!}-4\hat C\cdot\frac{\omega_h^m}{m!} \eeq
where $$\hat C=R_{ik\ell
j}\left(h^{\alpha\bar\beta}f^i_{\al}f^j_{\bbe}\right)
\left(h^{\gamma\bar\delta}f^k_{\ga}f^\ell_{\bar\delta}\right)$$ is
defined in (\ref{hatC}). From integration by parts, one obtains
 $$4\int_M |\p_\mathscr E\bp
f|^2\frac{\omega^m_h}{m!}-\int_M\hat C\cdot\frac{\omega_h^m}{m!}
=0.$$ If $(N,g)$ has  non-positive complex sectional curvature, then
$\hat C\leq 0$. Hence, we have $\hat C\equiv 0$ and $\p_\mathscr
E\bp f=0$, i.e. $f$ is pluri-harmonic.
 \eproof \elemma

\noindent\emph{The proof of Theorem \ref{main13}.} Let $(M,h)$ be a
compact astheno-K\"ahler manifold and $(N,g)$ be a Riemannian
manifold with non-positive complex sectional curvature. By Lemma
\ref{pluri}, every Hermitian harmonic map $f:(M,h)\>(N,g)$ is
pluri-harmonic. If $\sO_{T^*_M}(-1)$ is RC-positive, then by Theorem
\ref{main11}, the pluri-harmonic map $f:M\>(N,g)$ is constant. \qed

\vskip 2\baselineskip

\section{Further applications of the generalized energy
density}\label{further}

\noindent In this section, we discuss briefly some further
applications of the generalized energy, which are analogous to the
classical theory of harmonic maps. For more related topics, we refer
to the survey papers and books \cite{EL78, EL83, EL88, Xin96} and
the references therein.

\subsection{ RC-positivity for Riemannian curvature
tensor} As analogous to the  RC-positivity for abstract vector
bundles, one can define it for Riemannian manifolds. \bdefinition
Let $(N,g)$ be a Riemannian manifold. The curvature tensor $R$ of
$(N,g)$ is said to be \emph{RC-positive} if at each point $p\in N$
and for any nonzero vector $Z\in T_pN$, there exists a vector $W\in
T_pN$ such that \beq R(Z, W, W, Z)> 0. \eeq\edefinition \noindent
This terminology is a generalization of positive Riemannian
sectional curvature.  For instances, Riemannian manifolds with
positive Ricci curvature must be RC-positive. Similarly, one can
define the uniform RC-positivity and other similar notions.
\bdefinition Let $(N,g)$ be a Riemannian manifold. The curvature
tensor $R$ of $(N,g)$ is called \emph{uniformly RC-positive} if at
each point $p\in N$  there exists a vector $W\in T_pN$ such that for
any nonzero vector $Z\in T_pN$, \beq R(Z, W, W, Z)> 0.
\eeq\edefinition

\noindent We can define the energy density function $\mathscr Y$ on
the projective bundle $\P(T_M)\>M$
$$\mathscr Y=g_{ij}f^i_\alpha f^j_\beta \frac{W^\alpha W^\beta}{\sum h_{\gamma\delta} W^\gamma W^\delta}$$
 for a smooth map
$f:(M,h)\>(N,g)$ between Riemannian manifolds. By using this
setting, we can obtain similar  Hessian estimates as in Theorem
\ref{main0}, Proposition \ref{A1} and Proposition \ref{A2} for
totally geodesic maps (or some other harmonic maps), and rigidity of
such maps follow in a similar way.

\subsection{The extension of Yau's  function theory on
complete manifolds} It is easy to see that by using (\ref{S1}) and
(\ref{S-1}), their traces or integration by parts, we can extend
Yau's function theory (e.g. \cite{Yau75, Yau78}) on complete
manifolds by various generalized maximum principles. One of the key
steps in established in \cite[Corollary~2.3]{Yang17}.

\subsection{The first and second variations of the generalized
energy function} Let $f_t:(M,g)\>(N,h)$ be a family of smooth maps
 parameterized by $t\in\Delta$, the corresponding generalized energy
 density is denoted by $\mathscr Y_t$. The first and second
 variations of $\mathscr Y_t$ are powerful in analyzing the
 stability and related properties of harmonic maps as shown in the classical works
 \cite{SU81, SY80} and also a recent work \cite{FLW17}. The energy
 density (\ref{density2}) would be crucial in this context.

\subsection{The analytical extension of this method to hyperbolic
manifolds} By using the theory of RC-positivity and results in
algebraic geometry (in particular, seminal works of
Graber-Harris-Starr \cite{GHS03}, Boucksom-Demailly-Paun-Peternell
\cite{BDPP13} and Campana-Demailly-Peternell \cite{CDP14}), we
obtain the following rigidity theorem in
\cite[Corollary~1.5]{Yang18b}

\btheorem Let $(M,h)$ be a compact K\"ahler manifold with positive
holomorphic sectional curvature (or more generally, with uniformly
RC-positive tangent bundle). Then there is no non-constant
homomorphic map from $M$ to a Brody hyperbolic complex manifold $N$.
\etheorem

\noindent It is a natural task to explore a purely differential
geometric proof of this result. More generally, we propose the
following conjecture (see also similar versions in
\cite[Conjecture~1.9]{Yang18b}).

\noindent
\begin{conjecture}\label{con} Let $M$ and $N$ be a two compact complex
manifolds. If $\sO_{T^*M}(-1)$ is RC-positive and $N$ is Kobayashi
hyperbolic, then there is no non-constant holomorphic map from $M$
to $N$.
\end{conjecture}

\noindent
 The key difficulty is that
there is no Hermitian  metric with desire curvature positivity on
hyperbolic manifolds. It is a reasonable way to attack Conjecture
\ref{con} by using Proposition \ref{A2} and the Demailly-Semple
tower method (e.g. \cite{Dem18, BD18}).

\subsection{The generalized energy density on Grassmannian manifolds
$\mathrm{Gr}(k,T_M)$.} Let's  called that the curvature matrix is
call RC-positive if it has at least one positive eigenvalue. In this
case, we considered the projection of this positive direction on the
projective bundle $\P(T_M)$. If the curvature matrix has
$k$-positive directions (\cite{Yang18, Yang18b}),  we can consider
the associated Grassmannian manifold $\mathrm{Gr}(k,T_M)$. Many
results of this paper still work in this general setting.

\end{document}